\documentclass{article}
\usepackage{amsmath, amssymb, graphics}

\usepackage{amsthm}
\usepackage{amscd}

\newtheorem{thm}{Theorem}[section]
\newtheorem{df}{Definition}[section]
\newtheorem{rem}{Remark}[section]

\newtheorem{prp}{Proposition}[section]

\makeatletter
    
   \@addtoreset{figure}{chapter}
  \makeatother
  
\newcommand{\mathsym}[1]{{}}
\newcommand{\unicode}[1]{{}}

\makeatletter
    
    \@addtoreset{equation}{section}
  \makeatother
  
\begin{document}
\title{Fixed-Point Problems in Discrete Tomography: Case of Square Windows}
\author{Fumio HAZAMA\\Tokyo Denki University\\Hatoyama, Hiki-Gun, Saitama JAPAN\\
e-mail address:hazama@mail.dendai.ac.jp\\Phone number: (81)49-296-2911}
\date{\today}
\maketitle
\thispagestyle{empty}

\begin{abstract}
A kind of fixed-point problem in the area of discrete tomography is proposed and investigated. Our chief concern in this paper is the case of square windows in the plane. Dealing with the arrays which are bounded, of polynomial growth, and finite-ring-valued, one comes across several interesting phenomena of combinatorial and arithmetic nature.\\
keywords: discrete tomography; fixed point; square window; balanced array
\end{abstract}

\section{Introduction}
In the articles [2], [3], we show that various problems in the area of discrete tomography, for example as in [4], can be understood through the theory of distributions in a unified way. The main purpose of the present paper is to investigate several combinatorial phenomena which show up when we try to apply the general theory developed there to concrete problems. In particular we focus on the case of {\it square windows} and investigate what kind of arrays arise as solutions to tomographic problems. The shapes of solutions depend heavily on the conditions which we impose on the arrays. Accordingly we divide our presentation into several sections which deal with arrays that are {\it bounded}, {\it of polynomial growth}, $\mathbb{F}_p$-{\it valued}, and {\it balanced}. \\\\
The plan of the paper is as follows. In Section one, after we fix some notation and recall the main results of [2] and [3], we formulate the problem which is of our major concern throughout the paper. Section two is devoted to the study of the {\it characteristic polynomials} attached to the square and related windows. As a result we determine the structure of the solution spaces of arrays which are bounded. In Section three we consider our problem for the arrays of polynomial growth employing the tools developed in [3]. Section four deals with the case of $\mathbb{F}_p$-valued arrays, and Section five concerns with the case of $\mathbb{Z}_{n^2}$-valued balanced arrays. Here we come across to a problem which resembles to the "Sudoku" game. In this final section we propose several open problems to the reader.
\section{Problem Setting.}
For any commutative ring $R$, let $\mathbf{A}(R)=(R)^{\mathbb{Z}^2}$ denote the set of $R$-valued functions on $\mathbb{Z}^2$. We write elements of $\mathbf{A}(R)$ as $\mathbf{a}=(\mathbf{a}_{\mathbf{i}})$ with $\mathbf{i}=(i_1,i_2)\in\mathbb{Z}^2$, and call them $R$-{\it valued arrays} on $\mathbb{Z}^2$. For any non-empty finite subset $W$ of $\mathbb{Z}^2$ and for any $R$-valued array $\mathbf{a}\in\mathbb{A}(R)$, let $d_W(\mathbf{a})=\sum_{\mathbf{i}\in W}\mathbf{a}_{\mathbf{i}}\in R$, and call it {\it the degree of} $\mathbf{a}$ {\it with respect to} $W$. For any $\mathbf{k}\in\mathbb{Z}^2$, let $W+\mathbf{k}=\{\mathbf{i}+\mathbf{k};\mathbf{i}\in W\}$ denote the translation of $W$ by $\mathbf{k}$. By getting $d_{W+\mathbf{k}}(\mathbf{a}) (\mathbf{k}\in\mathbb{Z}^2)$ together, we obtain a new $R$-valued array $(d_{W+\mathbf{k}}(\mathbf{a}))_{\mathbf{k}\in\mathbb{Z}^2}$. We denote this new array by $\Delta_W(\mathbf{a})$ so that the set of $R$-valued arrays are equipped with a self-map $\Delta_W:\mathbf{A}(R)\rightarrow\mathbf{A}(R)$, which is easily seen to be an $R$-endomorphism through the natural $R$-module structure on $\mathbf{A}(R)$. The main purpose of this paper is to investigate the set of fixed points of $\Delta_W$ for various rings and finite subsets $W\subset\mathbb{Z}^2$. Accordingly we introduce the following notation:

\begin{df}
For any non-empty finite subset $W$ of $\mathbb{Z}^2$, let $Fix_W(R)$ denote the set of fixed points of $\Delta_W:\mathbf{A}(R)\rightarrow\mathbf{A}(R)$, namely
\begin{eqnarray*}
Fix_W(R)=\{\mathbf{a}\in \mathbf{A}(R);\Delta_W(\mathbf{a})=\mathbf{a}\}.
\end{eqnarray*}
\end{df}

\noindent
By the very definition of $\Delta_W$, we have the equality $Fix_W(R)=Fix_{W+\mathbf{k}}(R)$ for any $\mathbf{k}\in\mathbb{Z}^2$. Therefore we assume throughout the paper that $W$ contains the origin $O=(0,0)\in\mathbb{Z}^2$, and let $W^*=W\setminus \{O\}$. The following simple observation provides us with a connection between our problem and the results obtained in [2] and [3]:

\begin{prp}
For any non-empty finite subset $W$ of $\mathbb{Z}^2$, we have
\begin{eqnarray}
Fix_W(R)=\{\mathbf{a}\in\mathbf{A}(R);\Delta_{W^*}(\mathbf{a})=\mathbf{0}\},
\end{eqnarray}
where $\mathbf{0}$ denotes the all-zero array.
\end{prp}

\noindent
We will denote the set on the right hand side of (1,1) by $\mathbf{A}_{W^*}$. Here we recall one of the main results in [2]. For any window $W$, we put
\begin{eqnarray*}
m_W(x,y)=\sum_{(i_1,i_2)\in W}x^{i_1}y^{i_2}\in \mathbf{Z}[x^{\pm 1}, y^{\pm 1}],
\end{eqnarray*}
and call it the {\it characteristic (Laurent) polynomial} of $W$. For any Laurent polynomial $f$, we define $f^*$ by the following rule:
\begin{eqnarray*}
f^*(x,y)=f(1/x, 1/y).
\end{eqnarray*}
Furthermore we define a subset $\mathbf{A}^0_W(\mathbb{C})$ of $\mathbf{A}_W(\mathbb{C})$ by
\begin{eqnarray*}
&&\mathbf{A}^0_W(\mathbb{C})\\
&&=\{\mathbf{a}\in\mathbf{A}_W(\mathbb{C});\mbox{\it there exists a constant}\hspace{1mm}C\hspace{1mm}\mbox{\it such that }|\mathbf{a}_{\mathbf{i}}|<C(\mathbf{i}\in\mathbb{Z}^2)\}.
\end{eqnarray*}
Then we have the following equality [2, Theorem 3.2]
\begin{eqnarray}
\dim_{\mathbb{C}}\mathbf{A}^0_W(\mathbb{C})=\#(V_{\mathbf{T}^2}(m^*_W)),
\end{eqnarray}
where $\mathbf{T}=\{z\in\mathbb{C};|z|=1\}$ and $V_{\mathbf{T}^2}(m^*_W)$ denotes the zero locus of $m^*_W$ on $\mathbb{T}^2$.

\section{Square window}
For any integer $n\geq 2$, let $S(n)$ denote the subset of $\mathbb{Z}^2$ defined as $S(n)=\{(i_1,i_2)\in\mathbb{Z}^2;0\leq i_1, i_2\leq n-1\}$, and call it {\it the square window of width} $n$. Its characteristic polynomial $m_{S(n)}$ is given by
\begin{eqnarray*}
m_{S(n)}(x,y)=(1+x+x^2+\cdots +x^{n-1})(1+y+y^2+\cdots +y^{n-1}),
\end{eqnarray*}
hence that of $S(n)^*$ is given by
\begin{eqnarray*}
m_{S(n)^*}(x,y)=(1+x+x^2+\cdots +x^{n-1})(1+y+y^2+\cdots +y^{n-1})-1.
\end{eqnarray*}
For any positive integer $N$, we denote by $\mu_N$ the set of $N$-th roots of unity in $\mathbb{C}$, and let $\mu_N^*=\mu_N\setminus\{1\}$. 

\begin{prp}
For any $c\in\mu_{n-1}$ let
\begin{eqnarray*}
S_1(c)&=&\{(x,y)\in(\mu_{n-1}^*)^2;xy=c\},\\
S_2(c)&=&\{(x,y)\in(\mu_{(n-1)(n+1)}^*)^2;xy=c, x^{n+1}=c\}.
\end{eqnarray*}
Then we have
\begin{eqnarray*}
V_{\mathbf{T}^2}(m_{S(n)^*})=\bigcup_{c\in\mu_{n-1}}\left(S_1(c)\cup S_2(c)\right).
\end{eqnarray*}
\end{prp}

\noindent
{\it Proof}. Since $m_{S(n)^*}$ is defined over $\mathbb{R}$, if $(x_0,y_0)\in V_{\mathbf{T}^2}(m_{S(n)^*})$ then its complex conjugate $(\overline{x_0},\overline{y_0})$ is also in $V_{\mathbf{T}^2}(m_{S(n)^*})$. Furthermore since we have $\bar{z}=1/z$ for any $z\in\mathbf{T}$, we see that
\begin{eqnarray*}
V_{\mathbf{T}^2}(m_{S(n)^*})=V_{\mathbf{T}^2}(m_{S(n)^*})\cap V_{\mathbf{T}^2}(m^*_{S(n)^*}).
\end{eqnarray*}
Therefore we are to solve the following simultaneous equations:
\begin{eqnarray}
(1+x+x^2+\cdots +x^{n-1})(1+y+y^2+\cdots +y^{n-1})-1=0,\\
(1+x^{-1}+x^{-2}+\cdots +x^{-(n-1)})(1+y^{-1}+y^{-2}+\cdots +y^{-(n-1)})-1=0.
\end{eqnarray}
By multiplying (2.2) by $x^{n-1}y^{n-1}$ and subtracting it from (2.1), we obtain the equality
\begin{eqnarray*}
x^{n-1}y^{n-1}=1.
\end{eqnarray*} 
Hence there exists an element $c\in\mu_{n-1}$ such that $y=c/x$. First we deal with the case when $c\neq 1$. By inserting $y=c/x$ into the equation (2.1) and expanding it, we can compute the left hand side as follows:
\begin{eqnarray}
&&(1+x+x^2+\cdots +x^{n-1})\nonumber\\
&&+c(x^{-1}+1+x+\cdots +x^{n-2})\nonumber\\
&&+c^2(x^{-2}+x^{-1}+1+\cdots +x^{n-3})\nonumber\\
&&\hspace{5mm}\cdots\cdots\cdots\nonumber\\
&&+c^{n-1}(x^{-(n-1)}+x^{-(n-2)}+x^{-(n-3)}+\cdots +1)-1\nonumber\\
&&=x^{n-1}+\sum_{k=0}^1c^kx^{n-2}+\cdots+\left(\sum_{k=0}^{n-2}c^kx\right)+\left(\sum_{k=0}^{n-1}c^k-1\right)\nonumber\\
&&+\left(\sum_{k=1}^{n-1}c^kx^{-1}\right)+\sum_{k=2}^{n-1}c^kx^{-2}+\cdots+\sum_{k=n-2}^{n-1}c^kx^{-(n-2)}+c^{n-1}x^{-(n-1)}\nonumber\\
&&\mbox{(note here that the contents of three parentheses vanish since}\sum_{k=0}^{n-2}c^k=0)\nonumber\\
&&=x^{n-1}+\sum_{k=0}^1c^kx^{n-2}+\sum_{k=0}^2c^kx^{n-3}+\cdots+\sum_{k=0}^{n-3}c^kx^2\nonumber\\
&&+\sum_{k=2}^{n-1}c^kx^{-2}+\sum_{k=3}^{n-1}c^kx^{-3}+\cdots+\sum_{k=n-2}^{n-1}c^kx^{-(n-2)}+c^{n-1}x^{-(n-1)}\nonumber\\
&&=x^{n-1}+\sum_{k=0}^1c^kx^{n-2}+\sum_{k=0}^2c^kx^{n-3}+\cdots+\sum_{k=0}^{n-3}c^kx^2\nonumber\\
&&-c(x^{-2}+\sum_{k=0}^1c^kx^{-3}+\sum_{k=0}^2c^kx^{-4}+\cdots+\sum_{k=0}^{n-3}c^kx^{-(n-1)})\nonumber\\
&&=(1-cx^{-(n+1)})(x^{n-1}+\sum_{k=0}^1c^kx^{n-2}+\sum_{k=0}^2c^kx^{n-3}+\cdots+\sum_{k=0}^{n-3}c^kx^2)\nonumber\\
&&=(1-cx^{-(n+1)})x^2(x^{n-3}+\sum_{k=0}^1c^kx^{n-4}+\sum_{k=0}^2c^kx^{n-5}+\cdots+\sum_{k=0}^{n-3}c^k).
\end{eqnarray}
Since we have
\begin{eqnarray*}
&&(x-c)(x^{n-3}+\sum_{k=0}^1c^kx^{n-4}+\sum_{k=0}^2c^kx^{n-5}+\cdots+\sum_{k=0}^{n-3}c^k)\\
&=&x^{n-2}+x^{n-3}+\cdots+x+1,
\end{eqnarray*}
we see that the rightmost side of the above equalities (2.3) becomes 
\begin{eqnarray*}
(1-cx^{-(n+1)})x^2(x^{n-2}+x^{n-3}+\cdots+x+1)/(x-c).
\end{eqnarray*}
Note that if $c=\zeta_{n-1}^a$ for some $a$ with $1\leq a\leq n-2$, we have
\begin{eqnarray*}
(x^{n-2}+x^{n-3}+\cdots+x+1)/(x-c)=\prod_{\substack{1\leq k\leq n-2\\k\neq a}}(x-\zeta_{n-1}^k)
\end{eqnarray*}
Therefore the set of solutions of (2.1) in case of $c\in \mu_{n-1}^*$ is given by
\begin{eqnarray*}
&&\{(x,y)\in(\mu_{n-1}^*)^2;xy=c, x\neq c\}\\
&&\cup\{(x,y)\in(\mu_{(n-1)(n+1)}^*)^2;xy=c, x^{n+1}=c\}.
\end{eqnarray*}
Here we note that if $(x,y)\in(\mu_{n-1}^*)^2$ and $xy=c$, then we must have $x\neq c$.
Hence the set of solutions are simplified to be expressed as
\begin{eqnarray}
\{(x,y)\in(\mu_{n-1}^*)^2;xy=c\}\cup\{(x,y)\in(\mu_{(n-1)(n+1)}^*)^2;xy=c, x^{n+1}=c\}.
\end{eqnarray}

\noindent
Next we deal with case when $c=1$ so that $y=1/x$. Suppose that $x=1$ and hence $y=1$. Then the left hand side of the equation (2.1) becomes greater than or equal to three, since $n\geq 2$. It follows that $x\neq 1$. By inserting $y=1/x$ into (2.1) and multiplying it by $x^{n-1}(1-x)^2$, we have
\begin{eqnarray*}
&&x^{n-1}(1-x)^2\\
&&\times \{(1+x+x^2+\cdots +x^{n-1})(1+1/x+1/x^2+\cdots +1/x^{n-1})-1\}\\
&=&(1-x^n)^2-x^{n-1}(1-x)^2\\
&=&(1-2x^n+x^{2n})-(x^{n-1}-2x^n+x^{n+1})\\
&=&1-x^{n-1}-x^{n+1}+x^{2n}\\
&=&(1-x^{n-1})(1-x^{n+1})\\
&=&0
\end{eqnarray*}
Hence the set of solutions of (2.1) in case of $c=1$ is given by
\begin{eqnarray}
\{(x,y)\in\mu_{n-1}^*;xy=1\}\cup\{(x,y)\in\mu_{n+1}^*;xy=1\}.
\end{eqnarray}
Furthermore by comparing the set of solutions (2.4) and (2.5), we see that (2.5) is obtained by setting $c=1$ in (2.4). Hence we complete the proof.\qed\\\\
Example 2.1. When $n=2$, Proposition 2.1 implies that
\begin{eqnarray*}
V_{\mathbf{T}^2}(m_{S(2)^*})&=&\cup_{c\in\mu_1}\left(S_1(c)\cup S_2(c)\right)\\
&=&S_1(1)\cup S_2(1)\\
&=&\phi\cup\{(x,y)\in(\mu_3^*)^2;xy=1, x^3=1\}\\
&=&\{(\zeta_3,\zeta_3^2), (\zeta_3^2, \zeta_3)\}.
\end{eqnarray*}
Therefore the equality (1.2) implies that $\dim_{\mathbb{C}}\mathbf{A}_{S(2)^*}^0(\mathbb{C})=2$, and it follows from [2] that a basis of $\mathbf{A}_{S(2)^*}^0(\mathbb{C})$ is given by
\begin{eqnarray*}
\mathbf{a}^1&=&(\zeta_3^i\zeta_3^{2j})_{(i,j)\in\mathbb{Z}^2}=(\zeta_3^{i+2j})_{(i,j)\in\mathbb{Z}^2},\\
\mathbf{a}^2&=&(\zeta_3^{2i}\zeta_3^j)_{(i,j)\in\mathbb{Z}^2}=(\zeta_3^{2i+j})_{(i,j)\in\mathbb{Z}^2}.
\end{eqnarray*}
Furthermore, if we put
\begin{eqnarray*}
\mathbf{b}^1&=&\mathbf{a}^1+\mathbf{a}^2,\\
\mathbf{b}^2&=&\frac{1}{\zeta_3-\zeta_3^2}(\mathbf{a}^1-\mathbf{a}^2),
\end{eqnarray*}
then $\mathbf{b}^1$ and $\mathbf{b}^2$ constitute a $\mathbb{Q}$-basis of $\mathbf{A}_{S(2)^*}^0(\mathbb{Q})$. Specifically the values $\mathbf{b}^1_{(i,j)}, \mathbf{b}^2_{(i,j)}$ of these arrays at $(i,j)\in\mathbb{Z}^2$ are given by
\[
\mathbf{b}^1_{(i,j)}=
\left\{
\begin{array}{ll}
2, & i-j\equiv 0\pmod{3},\\
-1, & i-j\not\equiv 0\pmod{3},
\end{array}
\right.
\]
and
\[
\mathbf{b}^2_{(i,j)}=
\left\{
\begin{array}{ll}
0, & i-j\equiv 0\pmod{3},\\
1, & i-j\equiv 1\pmod{3},\\
-1, & i-j\equiv 2\pmod{3}.\\
\end{array}
\right.
\]
These arrays are depicted as follows:\\\\
$\mathbf{b}^1$:
 
\includegraphics{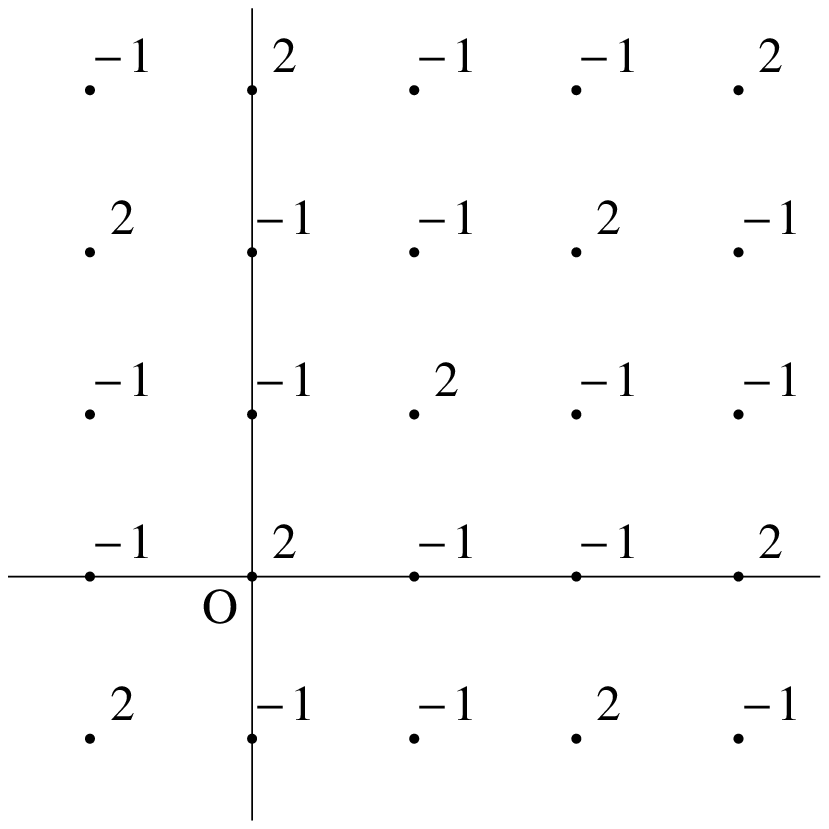}\\
$\mathbf{b}^2$:

\includegraphics{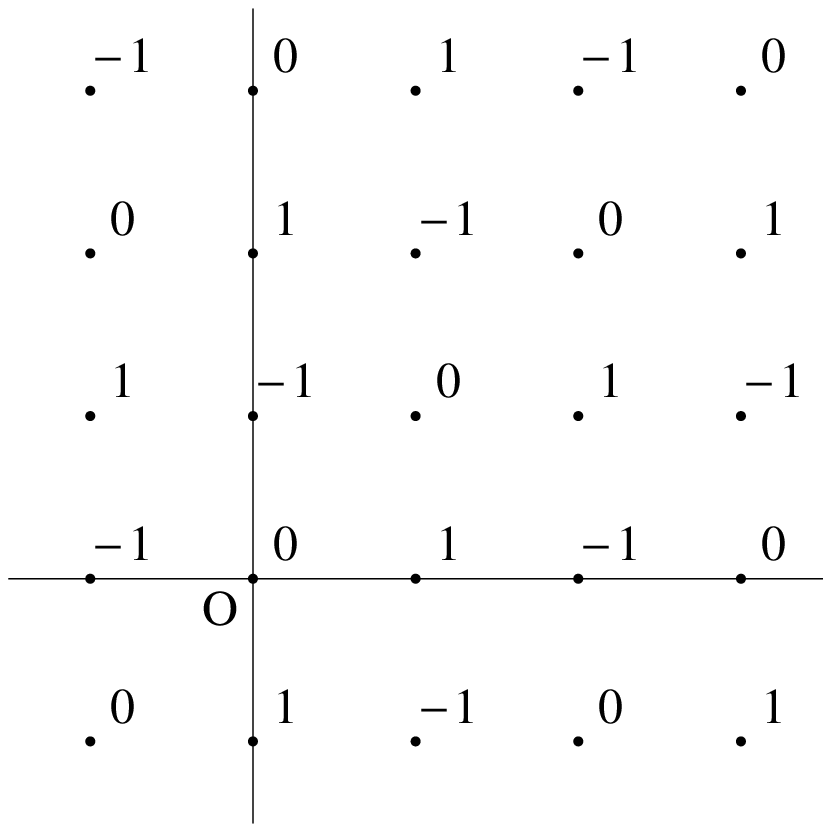}\\
We see from these figures that they are indeed fixed arrays for the square window $S(2)$.\\

\noindent
Example 2.2. When $n=3$, Proposition 2.1 implies that
\begin{eqnarray*}
V_{\mathbf{T}^2}(m_{S(3)^*})&=&\bigcup_{c\in\mu_2}\left(S_1(c)\cup S_2(c)\right)\\
&=&(S_1(1)\cup S_2(1))\cup(S_1(-1)\cup S_2(-1)),
\end{eqnarray*}
where
\begin{eqnarray*}
S_1(1)&=&\{(x,y)\in(\mu_2^*)^2;xy=1\}\\
&=&\{(-1,-1)\},\\
S_2(1)&=&\{(x,y)\in(\mu_8^*)^2;xy=1, x^4=1\}\\
&=&\{(i,-i), (-1,-1), (-i, i)\},\\
S_1(-1)&=&\{(x,y)\in(\mu_2^*)^2;xy=-1\}\\
&=&\phi,\\
S_2(-1)&=&\{(x,y)\in(\mu_8^*)^2;xy=-1, x^4=-1\}\\
&=&\{(\zeta_8, \zeta_8^3), (\zeta_8^3, \zeta_8), (\zeta_8^5, \zeta_8^7), (\zeta_8^7, \zeta_8^5)\}. 
\end{eqnarray*}
Hence we have
\begin{eqnarray*}
V_{\mathbf{T}^2}(m_{S(3)^*})=\{(-1,-1), (i,-i), (-i,i), (\zeta_8, \zeta_8^3), (\zeta_8^3, \zeta_8), (\zeta_8^5, \zeta_8^7), (\zeta_8^7, \zeta_8^5)\}. 
\end{eqnarray*}
We denote the arrays which correspond to the seven elements of $V_{\mathbf{T}^2}(m_{S(3)^*})$ by $\mathbf{a}^1, \mathbf{a}^2, \cdots, \mathbf{a}^7$, respectively. They constitute a basis of $\mathbf{A}_{S(3)^*}^0(\mathbb{C})$. Furthermore all of them are doubly-periodic, namely there exists a pair of linearly independent points $\mathbf{p}_k, \mathbf{q}_k\in\mathbb{Z}^2, 1\leq k\leq 7$, such that $\mathbf{a}^k_{\mathbf{i}+\mathbf{p}_k}=\mathbf{a}^k_{\mathbf{i}}$ and $\mathbf{a}^k_{\mathbf{i}+\mathbf{q}_k}=\mathbf{a}^k_{\mathbf{i}}$ hold.  More specifically, we see that
\begin{eqnarray*}
&&\mathbf{p}_1=(2,0), \mathbf{q}_1=(0,2),\\
&&\mathbf{p}_2=\mathbf{p}_3=(4,0), \mathbf{q}_2=\mathbf{q}_3=(0,4),\\
&&\mathbf{p}_4=\mathbf{p}_5=\mathbf{p}_6=\mathbf{p}_7=(8,0), \mathbf{q}_4=\mathbf{q}_5=\mathbf{q}_6=\mathbf{q}_7=(0,8).
\end{eqnarray*}
In order to make a $\mathbb{Q}$-basis, we let

\begin{eqnarray*}
\mathbf{b}^1&=&\mathbf{a}^1,\\
\mathbf{b}^2&=&\frac{1}{2}(\mathbf{a}^2+\mathbf{a}^3),\\
\mathbf{b}^3&=&\frac{1}{2i}(\mathbf{a}^2-\mathbf{a}^3),\\
\mathbf{b}^4&=&\frac{1}{4}(\mathbf{a}^4+\mathbf{a}^5+\mathbf{a}^6+\mathbf{a}^7),\\
\mathbf{b}^5&=&\frac{1}{4}(\zeta_8\mathbf{a}^4+\zeta_8^3\mathbf{a}^5+\zeta_8^5\mathbf{a}^6+\zeta_8^7\mathbf{a}^7),\\
\mathbf{b}^6&=&\frac{1}{4}(\zeta_8^2\mathbf{a}^4+\zeta_8^6\mathbf{a}^5+\zeta_8^2\mathbf{a}^6+\zeta_8^6\mathbf{a}^7),\\
\mathbf{b}^7&=&\frac{1}{4}(\zeta_8^3\mathbf{a}^4+\zeta_8\mathbf{a}^5+\zeta_8^7\mathbf{a}^6+\zeta_8^5\mathbf{a}^7).\\
\end{eqnarray*}

\noindent
Note that $\mathbf{b}^k$ inherit the same period from $\mathbf{a}^k$ for $1\leq k\leq 7$ by construction. Moreover if we define the operator of translation $T_{\mathbf{p}}, \mathbf{p}\in\mathbb{Z}^2$, of an array $\mathbf{a}$ by the rule $(T_{\mathbf{p}}(\mathbf{a}))_{\mathbf{i}}=\mathbf{a}_{\mathbf{i}-\mathbf{p}}$, then we see that
\begin{eqnarray*}
\mathbf{b}^3&=&T_{(1,0)}(\mathbf{b}^2),\\
\mathbf{b}^5&=&T_{(-1,0)}(\mathbf{b}^4),\\
\mathbf{b}^6&=&T_{(-2,0)}(\mathbf{b}^4),\\
\mathbf{b}^7&=&T_{(-3,0)}(\mathbf{b}^4).
\end{eqnarray*}
The arrays $\mathbf{b}^1, \mathbf{b}^2,$ and $ \mathbf{b}^4$ are depicted as follows:\\\\
$\mathbf{b}^1$:
 
\includegraphics{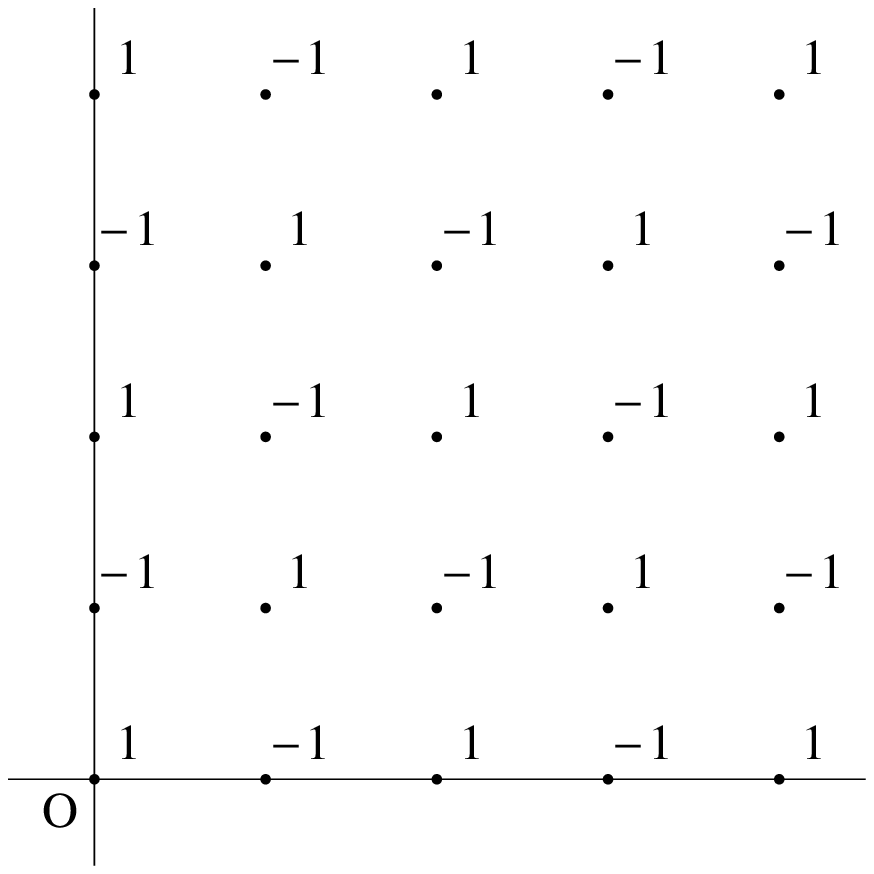}\\
$\mathbf{b}^2$:

\includegraphics{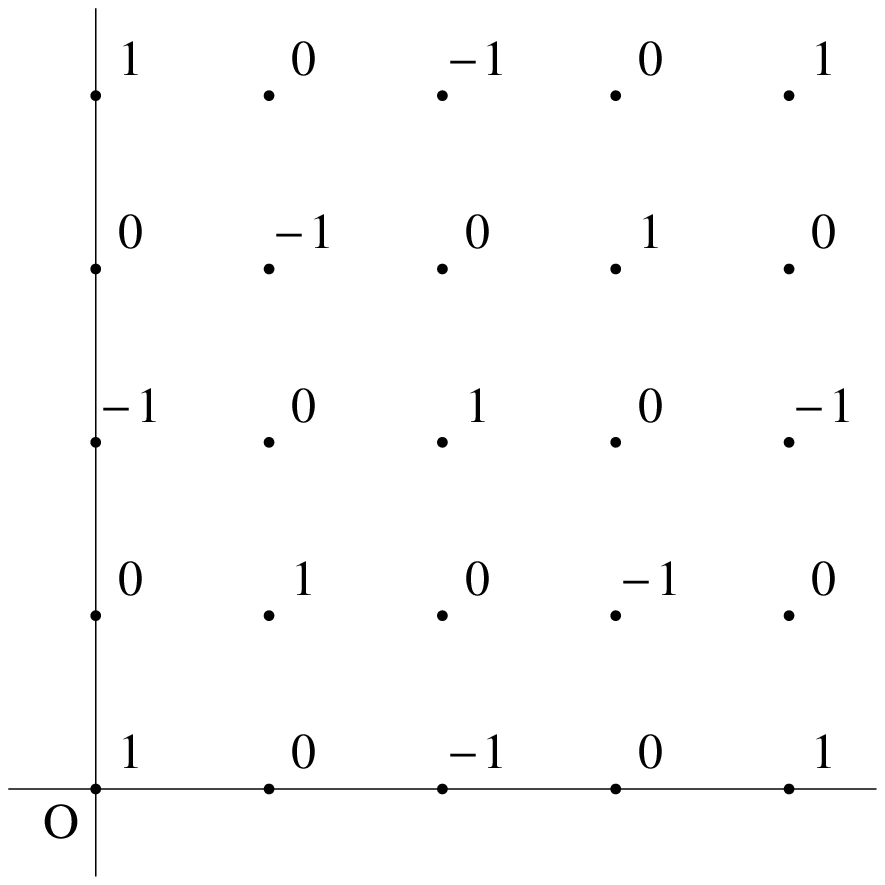}\\
$\mathbf{b}^4$:

\includegraphics{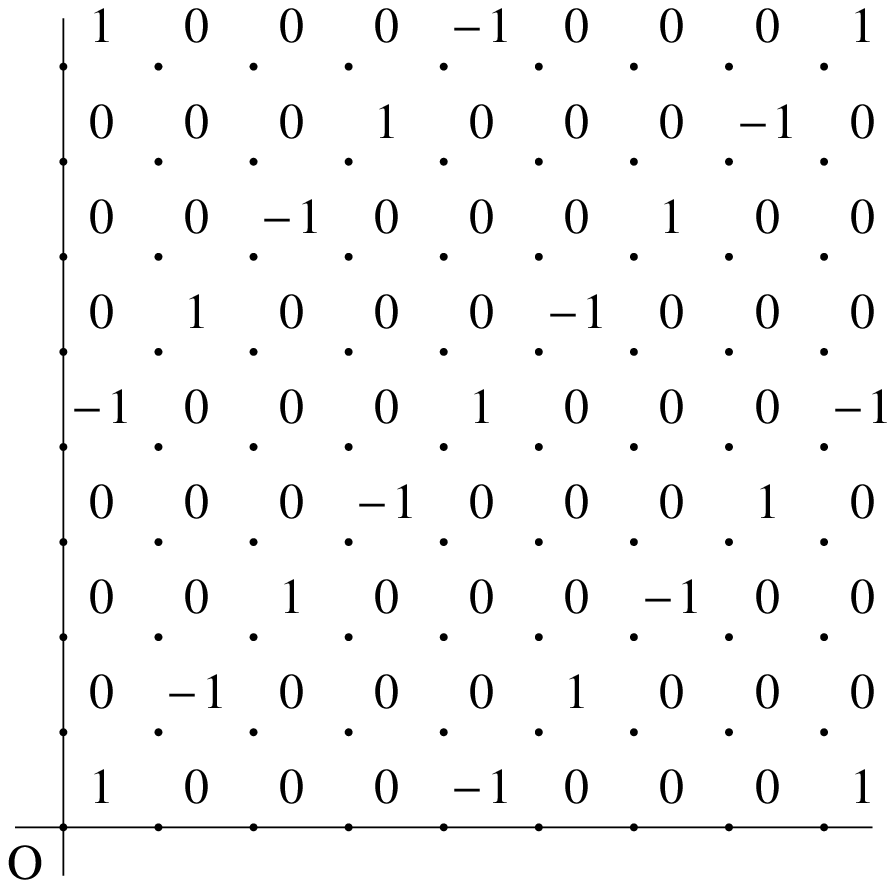}\\
We see from these figures that they are indeed fixed arrays for the square window $S(3)$.

\section{Unbounded arrays}
In this section we investigate unbounded arrays which are fixed by square windows. Here the theory developed in [3] plays a crucial role, and we recall some of main results of the paper for the convenience of the reader.\\

For any nonnegative integer $N$, let
\begin{eqnarray*}
\mathbf{A}(\mathbb{C})^N=\{\mathbf{a}\in\mathbf{A}(\mathbb{C});\mathbf{a}_{\mathbf{k}}=O(|\mathbf{k}|^N)\hspace{1mm}\mbox{as}\hspace{1mm}|\mathbf{k}|\rightarrow\infty\},
\end{eqnarray*}
and put
\begin{eqnarray*}
\mathbf{A}(\mathbb{C})^{poly}=\bigcup_{N\geq 0}\mathbf{A}(\mathbb{C})^N.
\end{eqnarray*}
Furthermore we let
\begin{eqnarray*}
\mathbf{A}(\mathbb{C})^N_W&=&\mathbf{A}_W\cap\mathbf{A}(\mathbb{C})^N,\\
\mathbf{A}(\mathbb{C})^{poly}_W&=&\mathbf{A}_W\cap\mathbf{A}(\mathbb{C})^{poly}.
\end{eqnarray*}
One of the main results in [3] relates the space $\mathbf{A}(\mathbb{C})^N_W$ with the space of solutions of a family of certain differential equations which are constructed from the data of $V_{\mathbb{T}^2}(m_W^*)$. In order to explain the connection, we need to introduce some notation. For each $\mathbf{p}\in V_{\mathbb{T}^2}(m_W^*)$, we fix a system of local coordinates $\mathbf{w}=(w_1, w_2)$ at $\mathbf{p}$. We put $P_2=\mathbb{C}[w_1, w_2]$, $D_2=\mathbb{C}[\partial_1, \partial_2]$, and let $W_2=\mathbb{C}\langle w_1,w_2,\partial_1, \partial_2\rangle$ denote the Weyl algebra in two variables, where $\partial_i=\partial/\partial w_i, 1\leq i\leq 2$. For any $D=\sum c_{\mathbf{i}}\partial^{\mathbf{i}}\in D_2$ with the usual convention for multi-index, let 
\begin{eqnarray*}
Sol(D)&=&\{g\in P_2;Dg=0\},\\
Sol(D)_{\leq N}&=&\{g\in P_2;Dg=0, \deg g\leq N\},\\
Sol(D)_{= N}&=&Sol(D)_{\leq N}/Sol(D)_{\leq N-1},
\end{eqnarray*} 
and
\begin{eqnarray*}
Supp(D)=\{\mathbf{i}\in\mathbb{Z}^2;c_{\mathbf{i}}\neq 0\}.
\end{eqnarray*} 
Let $F_-$ denote the $\mathbb{C}$-algebra automorphism of $W_2$ defined by the rule
\begin{eqnarray*}
F_-(w_i)&=&-\partial_i,\\
F_-(\partial_i)&=&w_i,
\end{eqnarray*}
for $i=1, 2$. Furthermore for any integers $k, s$, let $k^{\bar{s}}=(k+1)(k+2)\cdots (k+s)$, and for any $\mathbf{k}=(k_1, k_2), \mathbf{s}=(s_1, s_2)\in\mathbb{Z}^2$, let $\mathbf{k}^{\bar{\mathbf{s}}}=k_1^{\overline{s_1}}k_2^{\overline{s_2}}$. Then every array in $\mathbf{A}(\mathbb{C})^{poly}_W$ is obtained through the following five steps [3, Algorithm 3.3]:\\

\noindent
\hspace{3em}(A) {\it Find} $V_{\mathbb{T}^2}(m_W)$.\\
\hspace{3em}(B) {\it For} $\mathbf{p}\in V_{\mathbb{T}^2}(m_W)$, {\it let} $\mathbf{w}=\mathbf{z}-\mathbf{p}$ {\it and put} $m_W^{\mathbf{p}}=m_W(\mathbf{w}+\mathbf{p})$.\\
\hspace{3em}(C) {\it Put} $D_{\mathbf{p}}=F_-(m_W^{\mathbf{p}})$ {\it and determine} $Sol(D_{\mathbf{p}})$.\\
\hspace{3em}(D) {\it For any} $g\in Sol(D_{\mathbf{p}}$, {\it let} $F_-^{-1}(g)=\sum c_{\mathbf{s}}\partial^{\mathbf{s}}$.\\
\hspace{3em}(E) {\it Then} $\sum c_{\mathbf{s}}(-\mathbf{k})^{\bar{\mathbf{s}}}\mathbf{p}^{\mathbf{k}-\mathbf{s}}\in \mathbf{A}(\mathbb{C})^{poly}_W$.\\

\noindent
Furthermore we have the following dimension formula [3, Proposition 4.1] for the space of polynomial solutions. Here we use the notation $A-B$ for any pair of subsets of $\mathbb{Z}^2$ to mean the set $\{\mathbf{i}-\mathbf{j};(\mathbf{i}, \mathbf{j})\in A\times B\}$ of the differences of the elements of $A$ and $B$.
\begin{prp}
For any $D\in D_2$, we have
\begin{eqnarray*}
\dim Sol(D)_{\leq N}=\max\{\frac{(N+1)(N+2)}{2}-\#((\Delta_{\leq N}^2-supp(D))\cap\mathbb{Z}^2_{\geq 0}),0\}.
\end{eqnarray*}
\end{prp}

\noindent
We apply these to the square windows.\\

\noindent
Example 3.1. $\mathbf{A}(\mathbb{C})^{poly}_{S(2)^*}$: We have already examined the set which appears in step (A) in the previous section, and it follows from Example 2.1 that
\begin{eqnarray*}
V_{\mathbf{T}^2}(m_{S(2)^*})=\{(\zeta_3,\zeta_3^2), (\zeta_3^2, \zeta_3)\}.
\end{eqnarray*}
Let $\mathbf{p}^1=(\zeta_3,\zeta_3^2), \mathbf{p}^2=(\zeta_3^2, \zeta_3)$. Then for the polynomials in the step (B) we have
\begin{eqnarray*}
m_{S(2)^*}^{\mathbf{p}^1}&=&m_{S(2)^*}(\mathbf{w}+\mathbf{p}^1)\\
&=&(w_1+\zeta_3)+(w_1+\zeta_3)(w_2+\zeta_3^2)+(w_2+\zeta_3^2)\\
&=&-\zeta_3w_1-\zeta_3^2w_2+w_1w_2,\\
m_{S(2)^*}^{\mathbf{p}^2}&=&-\zeta_3^2w_1-\zeta_3w_2+w_1w_2.
\end{eqnarray*}
Hence the corresponding differential operators in the step (C) are given by
\begin{eqnarray*}
D_{\mathbf{p}^1}&=&\zeta_3\partial_1+\zeta_3^2\partial_2+\partial_1\partial_2,\\
D_{\mathbf{p}^2}&=&\zeta_3^2\partial_1+\zeta_3\partial_2+\partial_1\partial_2.\\
\end{eqnarray*}
Since both the supports of $D_{\mathbf{p}^1}$ and $D_{\mathbf{p}^2}$ are equal to $\{(1,0), (0,1), (1,1)\}$, it follows from Proposition 2.1 that
\begin{eqnarray*}
&&\dim Sol(D_{\mathbf{p}^1})_{\leq N}=\dim Sol(D_{\mathbf{p}^2})_{\leq N}\\
&&=\max\{\frac{(N+1)(N+2)}{2}-\#((\Delta_{\leq N}^2-\{(1,0), (0,1), (1,1)\})\cap\mathbb{Z}^2_{\geq 0}),0\}\\
&&=\max\{\frac{(N+1)(N+2)}{2}-\#(\Delta_{\leq N-1}^2),0\}\\
&&=\max\{\frac{(N+1)(N+2)}{2}-\frac{N(N+1)}{2},0\}\\
&&=N+1.
\end{eqnarray*}
Therefore we see that
\begin{eqnarray*}
\dim Sol(D_{\mathbf{p}^1})_{= N}=\dim Sol(D_{\mathbf{p}^2})_{= N}=1
\end{eqnarray*} 
for every $N\geq 0$. We illustrate below a representative of this quotient space for small $N$: (We use the symbols $x, y$ meaning $x=w_1, y=w_2$.)
\begin{table}[hbtp]
 \label{table:data_type}
 \begin{center}
  \begin{tabular}{lccc}
   \hline
   & $N = 0$  & $N = 1$  &  $N = 2$  \\
   \hline \hline
   $Sol(D_{\mathbf{p}^1})_{= N}$ & $const.$  & $(1)\hspace{1mm}\zeta_3^2x-\zeta_3y$ & $(2)\hspace{1mm}\zeta_3x^2-2xy+\zeta_3^2y^2-2x-2y$\\
   $Sol(D_{\mathbf{p}^2})_{= N}$ & $const.$  & $(3)\hspace{1mm}\zeta_3x-\zeta_3^2y$ & $(4)\hspace{1mm}\zeta_3^2x^2-2xy+\zeta_3y^2-2x-2y$\\
   \hline
  \end{tabular}
 \end{center}
\end{table}\\
By the step (E), each solution in this table gives rise to an array in $\mathbf{A}(\mathbb{C})^{poly}_{S(2)^*}$. Let $\mathbf{a}^{i}$ denote the array which corresponds to the solution $(i),\hspace{1mm}1\leq i\leq 4$, in the table above. (Note that the solution of degree zero provides us with the bounded arrays which are specified in the previous section.) Furthermore we put
\begin{eqnarray*}
\mathbf{b}^1&=&\mathbf{a}^1+\mathbf{a}^2,\\
\mathbf{b}^2&=&\mathbf{a}^1-\mathbf{a}^2,\\
\mathbf{b}^3&=&\frac{1}{2}(\mathbf{a}^3+\mathbf{a}^4),\\
\mathbf{b}^4&=&\frac{1}{6(2\zeta_3+1)}(\mathbf{a}^3-\mathbf{a}^4).\\
\end{eqnarray*} 

\noindent
These arrays are depicted as follows:\\\\
$\mathbf{b}^1$:
 
\includegraphics{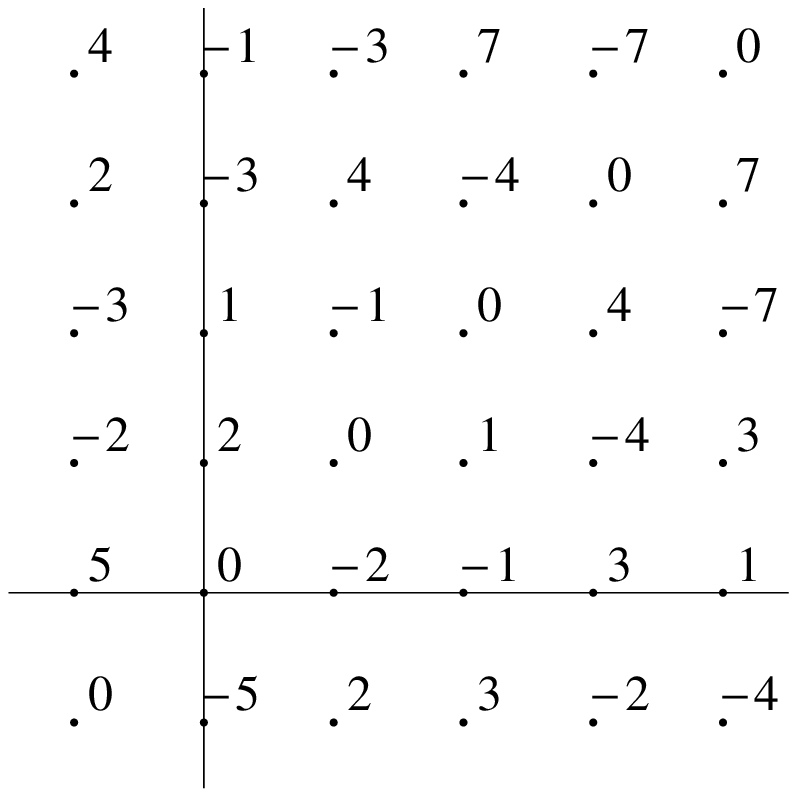}\\
$\mathbf{b}^2$:

\includegraphics{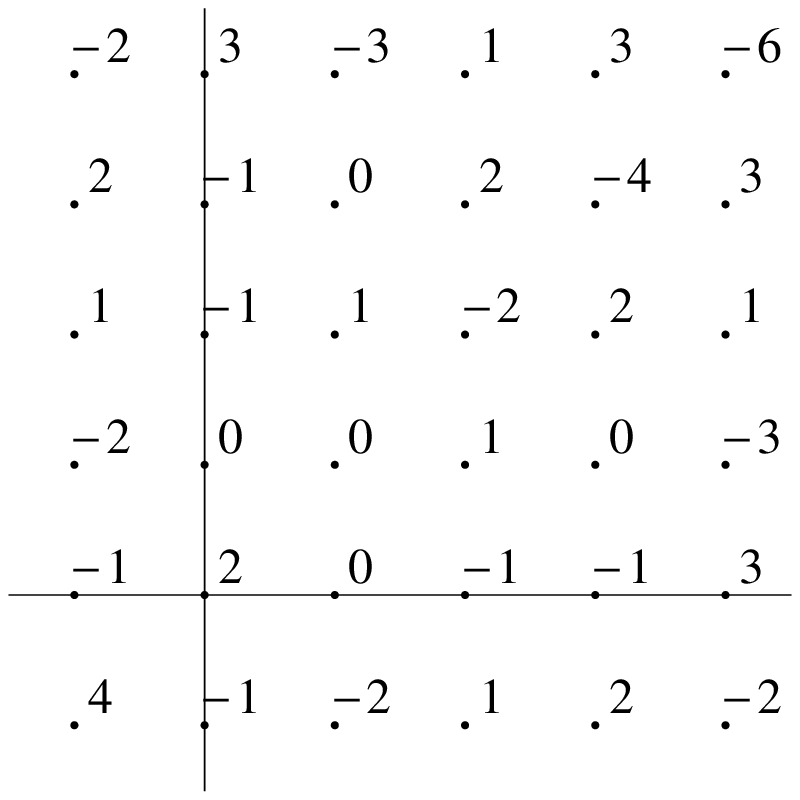}\\
$\mathbf{b}^3$:

\includegraphics{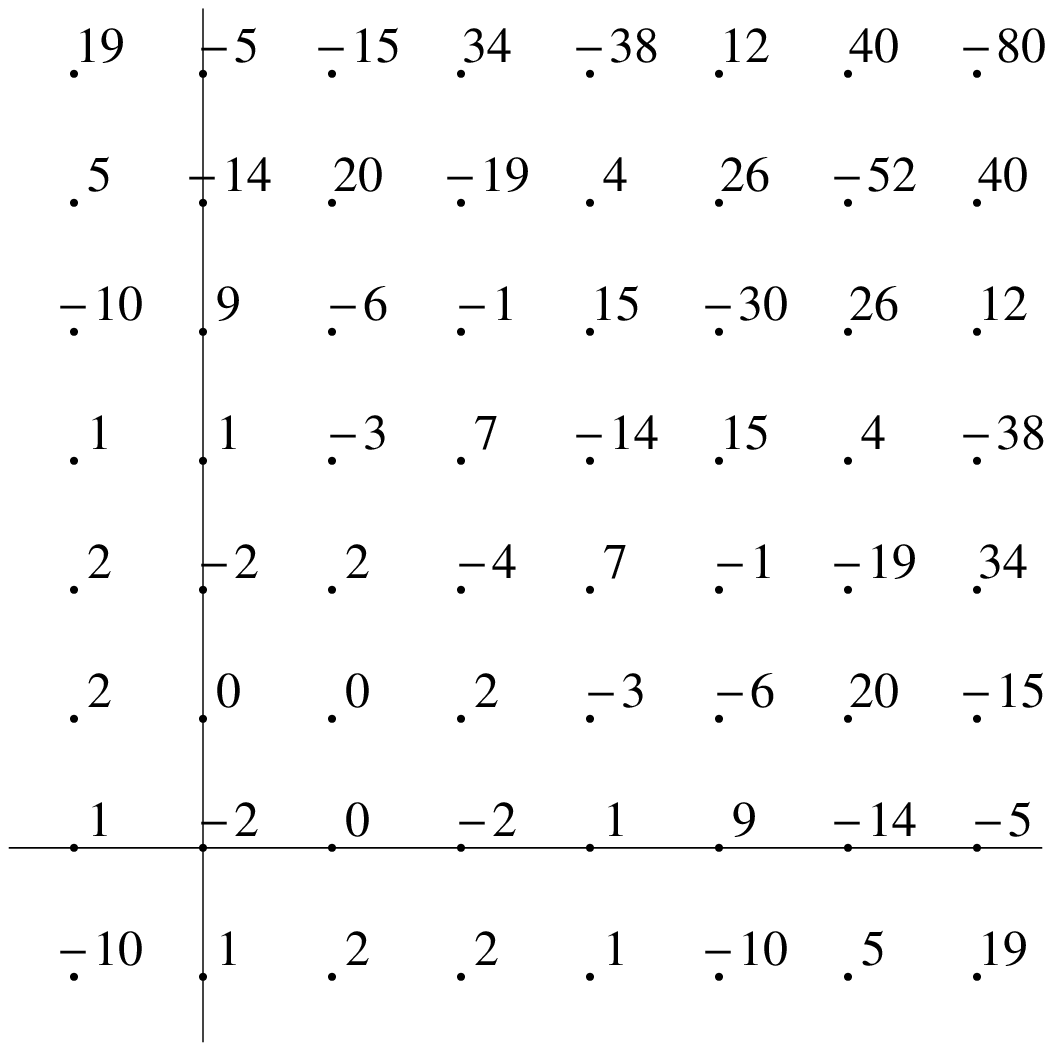}\\
$\mathbf{b}^4$:

\includegraphics{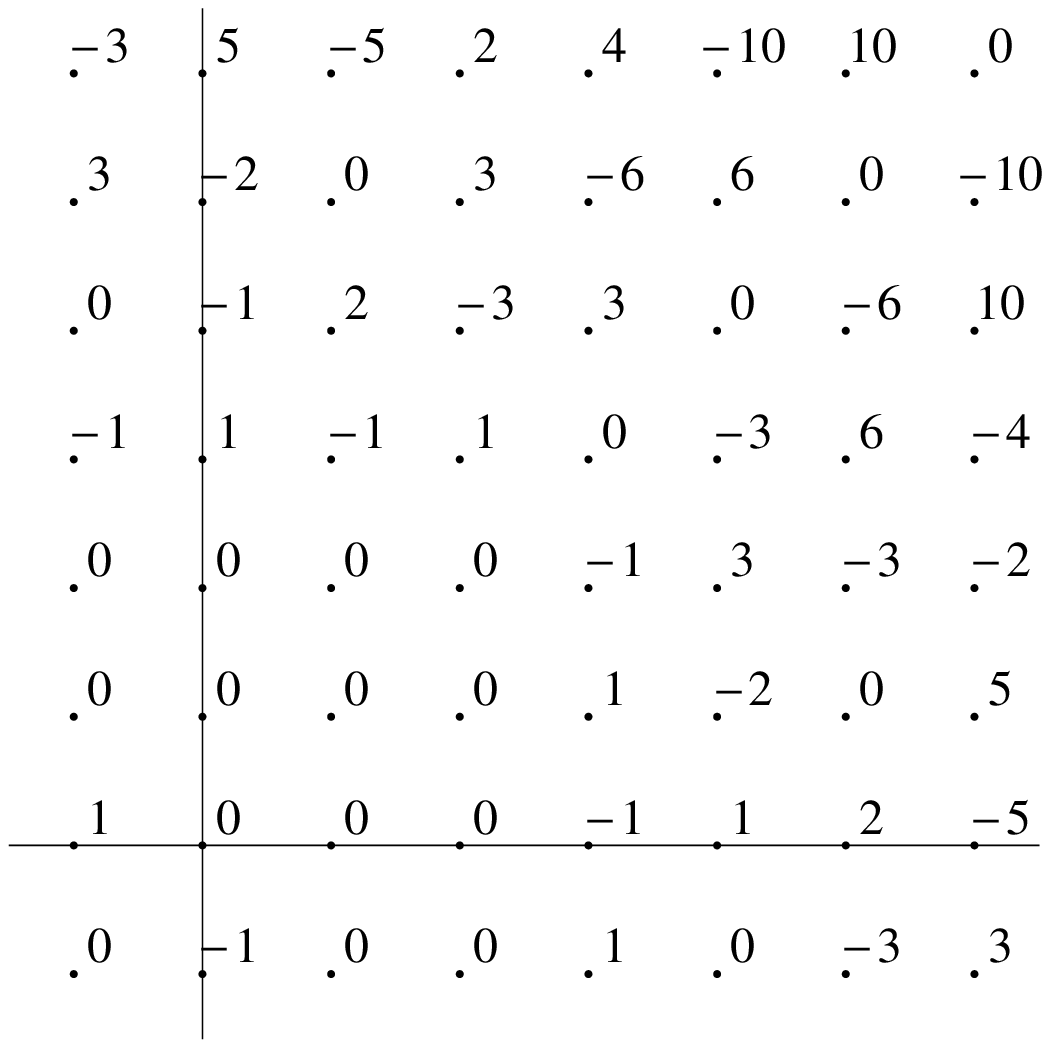}\\
We see directly from these figures that they are indeed fixed arrays for the square window $S(2)$.\\

\noindent
Example 3.2. $\mathbf{A}(\mathbb{C})^{poly}_{S(3)^*}$: It follows from Example 2.2 that
\begin{eqnarray*}
V_{\mathbf{T}^2}(m_{S(3)^*})=\{(-1,-1), (i,-i), (-i,i), (\zeta_8, \zeta_8^3), (\zeta_8^3, \zeta_8), (\zeta_8^5, \zeta_8^7), (\zeta_8^7, \zeta_8^5)\}. 
\end{eqnarray*}
We focus on the point $(-1, -1)\in V_{\mathbf{T}^2}(m_{S(3)^*})$, which gives rise to $\mathbb{Z}$-array in $\mathbf{A}(\mathbb{C})^{poly}_{S(3)^*}$. The differential equation in the Step (C) is found to be
\begin{eqnarray*}
D_{(-1,-1)}f=((1+\partial_1+\partial_1^2)(1+\partial_2+\partial_2^2)-1)f=0,
\end{eqnarray*}
and a representative of $Sol(D_{(-1,-1)})_{= 1}$ (resp. $Sol(D_{(-1,-1)})_{= 2}$) is given by $x-y$ (resp. $(x-y)^2-(x+y)$). The corresponding arrays are depicted as follows:\\\\
\noindent

\includegraphics{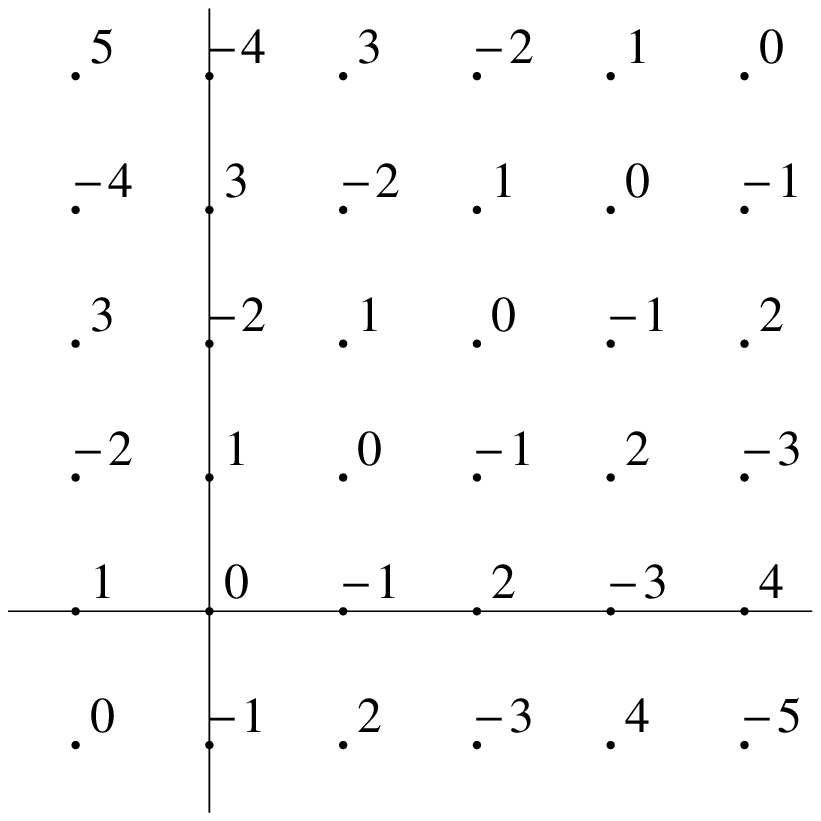}

\includegraphics{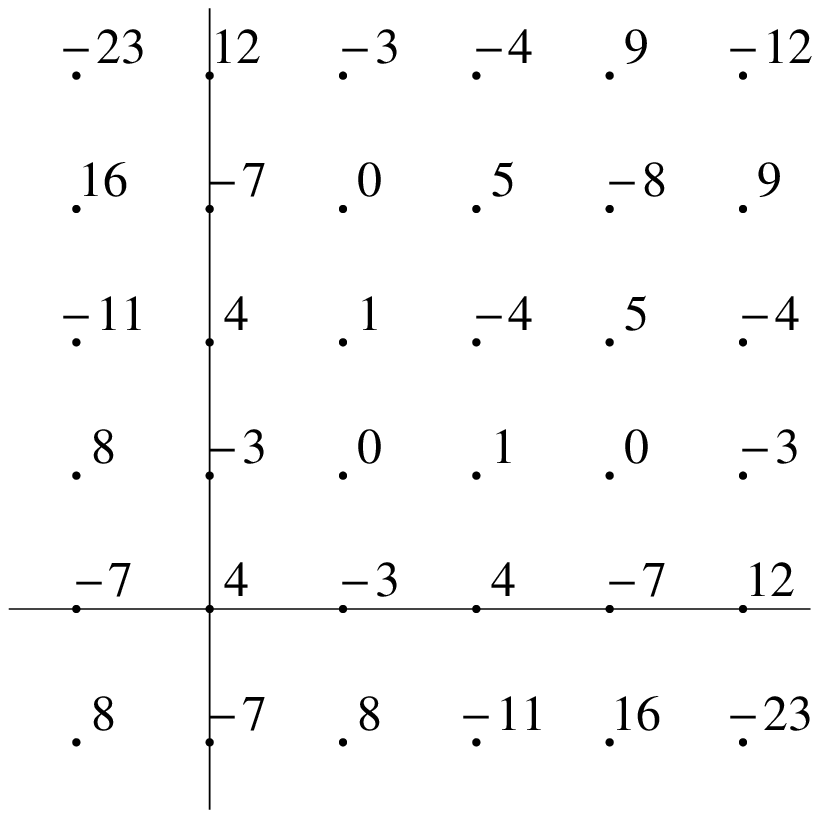}

\noindent
We can readily see that they are indeed fixed arrays for the square window $S(3)$.\\

\section{Square windows on the $p-$torus}
For any odd prime number $p$, let $T_p=\mathbb{F}_p\times\mathbb{F}_p$, and call it the $p$-{\it torus}. Let $\mathbf{A}_p(\mathbb{F}_p)=(\mathbb{F}_p)^{T_p}$ denote the set of $\mathbb{F}_p$-valued functions on the $p$-torus, and call its elements $p$-{\it arrays}. We can and will identify the set of $p$-arrays with the subset of $\mathbf{A}(\mathbb{F}_p)$ which have period $(p,0)$ and $(0,p)$. With this identification understood, the main object of our study in this section is the set $Fix_{S(n)}(\mathbb{F}_p)\cap\mathbf{A}_p(\mathbb{F}_p)$ for every $n$ with $1<n<p$. As is explained in Section one, this set is equal to $\mathbf{A}_{S(n)^*}\cap\mathbf{A}_p(\mathbb{F}_p)$, which will be denoted by $\mathbf{A}_p(\mathbb{F}_p)_{S(n)^*}$. In order to investigate the structure of this space we cannot apply the theory developed in [2], [3], since it deals with $\mathbb{C}$-valued functions through the theory of distributions. We employ here the theory of {\it group determinants}. As will be seen in the course of our description, the theory looks like a {\it finite-field-analogue} of that of distributions. We illustrate our method by the smallest possible example of $\mathbf{A}_3(\mathbb{F}_3)_{S(2)^*}$.\\

\noindent
Example 4.1. $\mathbf{A}_3(\mathbb{F}_3)_{S(2)^*}$. For any integer $k$ with $0\leq k\leq 8$, let $\mathbf{e}_k$ denote the element of $\mathbf{A}_3(\mathbb{F}_3)$ defined by
\[
\mathbf{e}_k(i,j)=
\left\{
\begin{array}{ll}
1, & k=i+3j,\\
0, & otherwise.
\end{array}
\right.
\]
Here we identify the elements of $\mathbf{F}_3$ with their least nonnegative representatives in $\mathbb{Z}$. With respect to the bases $(\mathbf{e}_k)_{0\leq k\leq 8}$, the representation matrix $M_{S(2)^*}$ of the linear operator $\Delta_{S(2)*}$ on  $\mathbf{A}_3(\mathbb{F}_3)$ is computed as
\begin{eqnarray*}
M_{S(2)^*}=\left(
\begin{array}{ccccccccc}
 0 & 1 & 0 & 1 & 1 & 0 & 0 & 0 & 0 \\
 0 & 0 & 1 & 0 & 1 & 1 & 0 & 0 & 0 \\
 1 & 0 & 0 & 1 & 0 & 1 & 0 & 0 & 0 \\
 0 & 0 & 0 & 0 & 1 & 0 & 1 & 1 & 0 \\
 0 & 0 & 0 & 0 & 0 & 1 & 0 & 1 & 1 \\
 0 & 0 & 0 & 1 & 0 & 0 & 1 & 0 & 1 \\
 1 & 1 & 0 & 0 & 0 & 0 & 0 & 1 & 0 \\
 0 & 1 & 1 & 0 & 0 & 0 & 0 & 0 & 1 \\
 1 & 0 & 1 & 0 & 0 & 0 & 1 & 0 & 0
\end{array}
\right).
\end{eqnarray*}
Its row-reduced form modulo 3 is found to be
\begin{eqnarray*}
\left(
\begin{array}{ccccccccc}
 1 & 0 & 0 & 0 & 0 & 0 & 2 & 2 & 1 \\
 0 & 1 & 0 & 0 & 0 & 0 & 1 & 2 & 2 \\
 0 & 0 & 1 & 0 & 0 & 0 & 2 & 1 & 2 \\
 0 & 0 & 0 & 1 & 0 & 0 & 1 & 0 & 1 \\
 0 & 0 & 0 & 0 & 1 & 0 & 1 & 1 & 0 \\
 0 & 0 & 0 & 0 & 0 & 1 & 0 & 1 & 1 \\
 0 & 0 & 0 & 0 & 0 & 0 & 0 & 0 & 0 \\
 0 & 0 & 0 & 0 & 0 & 0 & 0 & 0 & 0 \\
 0 & 0 & 0 & 0 & 0 & 0 & 0 & 0 & 0 \\
\end{array}
\right),
\end{eqnarray*}
and hence we have 
\begin{eqnarray*}
\dim_{\mathbb{F}_3}\mathbf{A}_3(\mathbb{F}_3)_{S(2)^*}=3.
\end{eqnarray*}
The corresponding three bases, which we call $\mathbf{a}^1, \mathbf{a}^2$, and $\mathbf{a}^3$,  are depicted as follows:\\\\
\noindent
$\mathbf{a}^1$:

\includegraphics{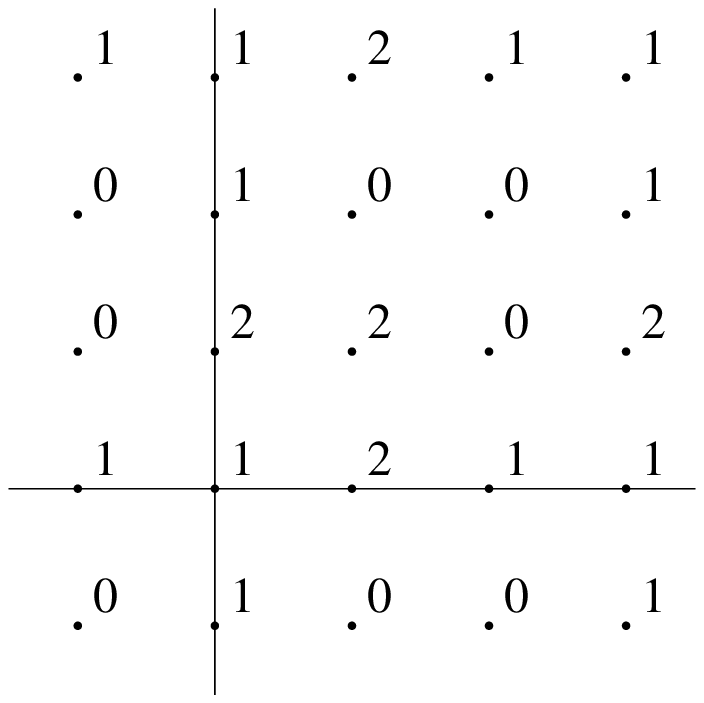}\\
$\mathbf{a}^2$:

\includegraphics{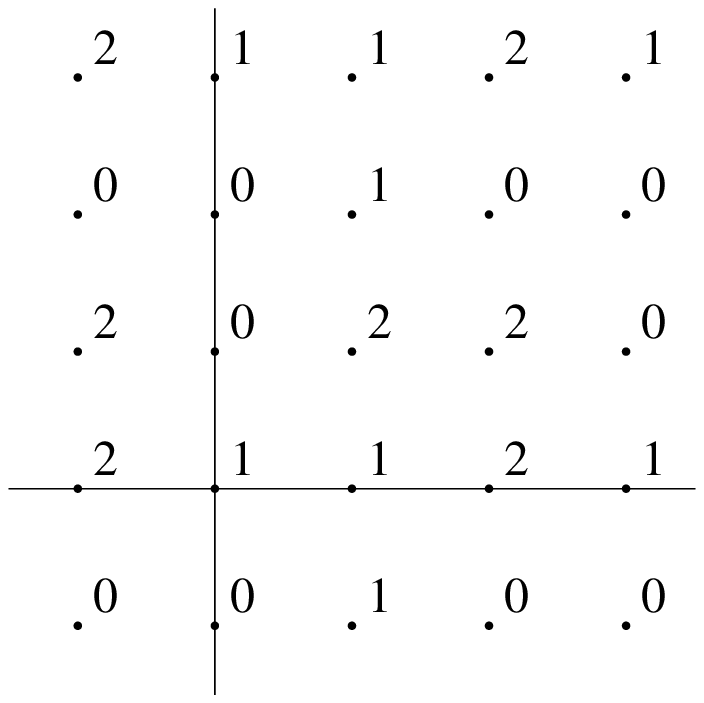}\\
$\mathbf{a}^3$:

\includegraphics{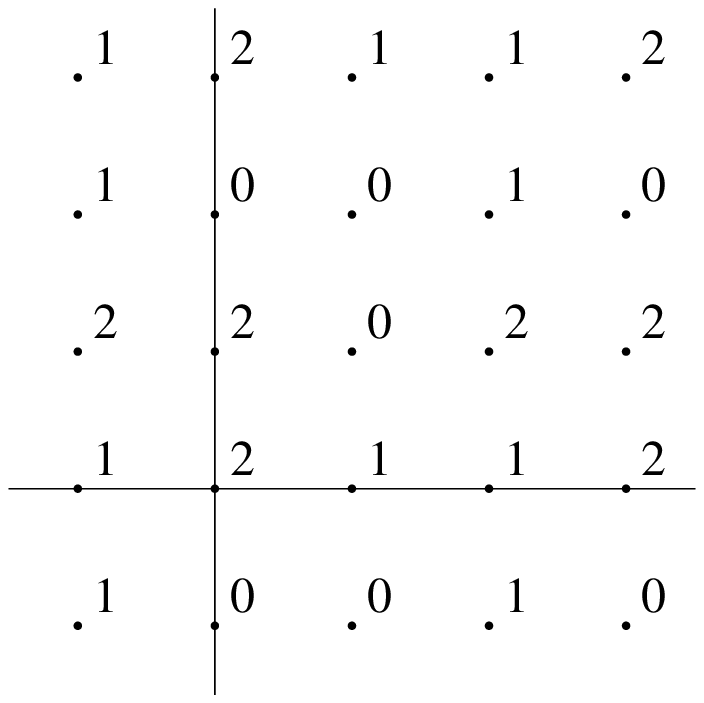}\\
As is evident from these figures, we have
\begin{eqnarray*}
\mathbf{a}^2&=&T_{(1,0)}\mathbf{a}^1,\\
\mathbf{a}^3&=&T_{(2,0)}\mathbf{a}^1.
\end{eqnarray*}
On the other hand, we have seen in Example 2.1 that 
\begin{eqnarray*}
\dim_{\mathbb{Q}}\mathbf{A}_{S(2)^*}^0(\mathbb{Q})=2,
\end{eqnarray*}
and every array in $\mathbf{A}_{S(n)^*}^0(\mathbb{Q})$ has periods $(3,0)$ and $(0,3)$. Hence by reducing modulo 3 the $\mathbb{Z}$-arrays of this space, we obtain the arrays in $\mathbf{A}_3(\mathbb{F}_3)_{S(2)^*}$. Indeed the bases $\mathbf{b}^1$ and $\mathbf{b}^2$ reduce to 
\begin{eqnarray*}
\mathbf{b}^1\pmod 3&=&2(\mathbf{a}^1+\mathbf{a}^2+\mathbf{a}^3),\\
\mathbf{b}^2\pmod 3&=&\mathbf{a}^1-\mathbf{a}^2.
\end{eqnarray*}
In contrast to this example, the case when $p>3$ is completely different:

\begin{prp}
For any prime $p>3$, we have
\begin{eqnarray*}
\mathbf{A}_p(\mathbb{F}_p)_{S(2)^*}=\{0\}.
\end{eqnarray*}
\end{prp}

\noindent
{\it Proof}. As is done in Example 4.1, for any integer $k$ with $0\leq k\leq p^2-1$, let $\mathbf{e}_k$ denote the element of $\mathbf{A}_p(\mathbb{F}_p)$ defined by
\[
\mathbf{e}_k(i,j)=
\left\{
\begin{array}{ll}
1, & k=i+pj,\\
0, & otherwise.
\end{array}
\right.
\]
These arrays contitute a basis of $\mathbf{A}_p(\mathbb{F}_p)$. With respect to this basis, we denote by $M_{S(2)^*}$ the representation matrix of the linear operator $\Delta_{S(2)*}$ on  $\mathbf{A}_p(\mathbb{F}_p)$. We will show that $\det M_{S(2)^*}\neq 0$, as an element of $\mathbb{F}_p$, by using the theory of group determinants [1]. For any finite group $G=\{g_1,\cdots,g_n\}$ with $g_1$ the identity element, the group matrix $M_G(x_{g_1},\cdots,x_{g_n})$ with variable $x_{g_k},\hspace{1mm}1\leq k\leq n$, is defined to be the $n\times n$-matrix whose $(i,j)$-entry is equal to $x_{g_ig_j^{-1}}$. For any $i\in\mathbb{F}_p$, let $\chi_i\in {\rm Hom}(\mathbb{F}_p,\mathbb{C}^*)$ denote the additive character of $\mathbb{F}_p$ defined by
\begin{eqnarray*}
\chi_i(a)=\zeta_p^{ai}.
\end{eqnarray*}
These characters constitute a complete set of characters of $\mathbb{F}_p$. One of main results of the theory of group determinants says that the determinant of the matrix $M=M_{T_p}(x_{(0,0)}, x_{(0,1)},\cdots,x_{(p-2,p-1)},x_{(p-1,p-1)})$ is factored as follows:
\begin{eqnarray}
\det M=\prod_{(i,j)\in T_p}\sum_{(k,\ell)\in T_p}\chi_i(k)\chi_j(\ell)x_{(k,\ell)}.
\end{eqnarray}
We can apply this fact to our matrix $M_{S(2)^*}$, since it coincides with the group matrix $M_{T_p}(x_{(0,0)}, x_{(0,1)},\cdots,x_{(p-2,p-1)},x_{(p-1,p-1)})$ with specialization
\[
x_{(i,j)}=
\left\{
\begin{array}{ll}
1, & -(i,j)\equiv (1,0), (0,1), (1,1)\pmod p,\\
0, & otherwise.
\end{array}
\right.
\]
Since the map $m_{-1}:T_p\rightarrow T_p$ defined by $(i,j)\mapsto -(i,j)\pmod p$ is a bijection on $T_p$ (actually an automorphism), the product on the right hand side in (4.1) can be expressed as
\begin{eqnarray*}
\prod_{(i,j)\in T_p}\sum_{(k,\ell)\in T_p}\chi_i(k)\chi_j(\ell)x_{(k,\ell)}=\prod_{(i,j)\in T_p}\sum_{(k,\ell)\in T_p}\chi_{p-i}(k)\chi_{p-j}(\ell)x_{(k,\ell)}.
\end{eqnarray*}
Hence we have
\begin{eqnarray*}
&&\det M_{S(2)^*}\\
&&=\prod_{(i,j)\in T_p}(\chi_i(1)\chi_j(0)x_{(1,0)}+\chi_i(0)\chi_j(1)x_{(0,1)}+\chi_i(1)\chi_j(1)x_{(1,1)})\\
&&=\prod_{(i,j)\in T_p}(\zeta_p^i+\zeta_p^j+\zeta_p^{i+j})
\end{eqnarray*}
Furthermore since the polynomial $X^p-1$ is factored as $(X-1)^p$ in $\mathbb{F}_p[X]$, the primitive $p$-th root of unity $\zeta_p$ gives rise to $1\pmod p$. Hence 
\begin{eqnarray*}
\det M_{S(2)^*}&\equiv& \prod_{(i,j)\in T_p}(1+1+1)\\
&\equiv&3^{p^2}\\
&\equiv&3\pmod p,
\end{eqnarray*}
the last line being a consequence of the Fermat's little theorem. This congruence means that the representation matrix $M_{S(2)^*}$ of the operator $\Delta_{S(2)^*}:T_p\rightarrow T_p$ is invertiable modulo $p$. This completes the proof. \qed\\\\

\noindent
Only a bit of change of the poof above enables one to obtain the most general result about $\mathbf{A}_p(\mathbb{F}_p)_{S(n)^*}$:

\begin{thm}
For any odd prime and for any integer $n$ with $2\leq n\leq p-2$, we have
 \begin{eqnarray*}
\mathbf{A}_p(\mathbb{F}_p)_{S(n)^*}=\{0\}.
\end{eqnarray*}
On the other hand when $n=p-1$, we have
\begin{eqnarray*}
\mathbf{A}_p(\mathbb{F}_p)_{S(n)^*}\neq\{0\}.
\end{eqnarray*}
\end{thm}

\noindent
{\it Proof}. For any $n$ with $2\leq n\leq p-1$, by a similar argument to the one given above, we see that
\begin{eqnarray*}
\det M_{S(n)^*}&\equiv& \prod_{(i,j)\in T_p}\#(S(n)^*)\\
&\equiv&(n^2-1)^{p^2}\\
&\equiv&(n^2-1)\\
&\equiv&(n-1)(n+1)\pmod p.
\end{eqnarray*}
Hence the matrix $M_{S(n)^*}$ is invertible if and only if $n+1=p$, namely $n=p-1$. This completes the proof. \qed\\\\

\section{Balanced arrays on the $n$-torus}
As is seen in Theorem 4.1, the square window $S(p-1)^*$ has always zero-sum array on the $p$-torus. Indeed the all-one array on $T_p$ gives rise to such an array, for the number of elements in $S(p-1)^*$ is equal to $(p-1)^2-1=p^2-2p$, which is equal to zero modulo $p$. In this section we investigate what occurs if we impose a strict condition on the entries of arrays.\\

For any integer $n\geq 3$, let $T_n=\mathbb{Z}_n\times\mathbb{Z}_n$ and call it the $n$-torus, and let $\mathbf{A}_n(\mathbb{Z}_{n^2})$ denote the set of $\mathbb{Z}_{n^2}$-valued functions on $T_n$. Let $\mathbf{A}_n(\mathbb{Z}_{n^2})_{S(n-1)^*}$ denote the subspace of $\mathbf{A}_n(\mathbb{Z}_{n^2})$ consisting of zero-sum arrays with respect to the window $S(n-1)^*$. We say an array $\mathbf{a}$ in $\mathbf{A}_n(\mathbb{Z}_{n^2})$ to be balanced if
\begin{eqnarray*}
\{\mathbf{a}_{(i,j)};(i,j)\in T_n\}=\mathbb{Z}_{n^2},
\end{eqnarray*}
namely if all the values of the array $\mathbf{a}$ are distinct. We consider the following problem:\\

\noindent
(P) {\it Does there exists a balanced array on} $T_n$ {\it such that it is a zero-sum array for the window} $S(k)^*,\hspace{1mm}2\leq k\leq n-1$, {\it too}?\\

\noindent
We illustrate the situation by the case when $n=3$ and $k=2$.\\

\noindent
Example 5.1. Window $S(2)^*$ on the 3-torus $T_3$: Let $\mathbf{x}$ denote the most general array on $T_3$ whose entries are defined  by $\mathbf{x}_{(i,j)}=x_{i+3j}$ for every $(i,j)\in T_3$. This is depicted below:

\includegraphics{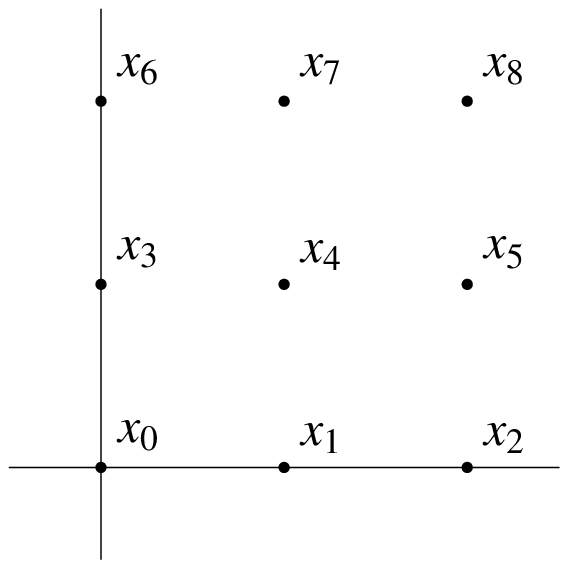}\\
We represent this array as the row vector $(x_0, x_1, \cdots, x_8)$. We find with the help of a computer that there are 12 balanced arrays in $\mathbf{A}_3(\mathbb{Z}_{3^2})_{S(2)^*}$ with $x_0=0$. We name them as follows;
\begin{eqnarray}
\mathbf{a}^1=(0,1,2,5,3,4,7,8,6), &&\mathbf{a}^2=(0,1,5,2,6,7,4,8,3), \nonumber\\
\mathbf{a}^3=(0,2,1,4,3,5,8,7,6), &&\mathbf{a}^4=(0,2,4,1,6,8,5,7,3),\nonumber\\
\mathbf{a}^5=(0,4,2,8,6,1,7,5,3), &&\mathbf{a}^6=(0,4,8,2,3,7,1,5,6),\nonumber\\
\mathbf{a}^7=(0,5,1,7,6,2,8,4,3), &&\mathbf{a}^8=(0,5,7,1,3,8,2,4,6),\nonumber\\
\mathbf{a}^9=(0,7,5,8,3,1,4,2,6), &&\mathbf{a}^{10}=(0,7,8,5,6,4,1,2,3),\nonumber\\
\mathbf{a}^{11}=(0,8,4,7,3,2,5,1,6), &&\mathbf{a}^{12}=(0,8,7,4,6,5,2,1,3).
\end{eqnarray}
We will show that this set is acted upon by a certain group. A permutation $\pi$ on the eight letters $\{1, 2, \cdots, 8\}$ acts on the set of arrays on $T_3$ through the rule: 
\begin{eqnarray*}
(x_0, x_1, \cdots, x_8)\mapsto (x_0, x_{\pi^{-1}(1)}, \cdots, x_{\pi^{-1}(8)}).
\end{eqnarray*}
Let us introduce a permutation $p$ defined by
\begin{eqnarray*}
p=\left(
\begin{array}{cccccccc}
 1 & 2 & 3 & 4 & 5 & 6 & 7 & 8 \\
 3 & 1 & 6 & 8 & 2 & 7 & 5 & 4
\end{array}
\right)
\end{eqnarray*}
Then the array $\mathbf{x}$ is transformed through $p$ to the following array:

\includegraphics{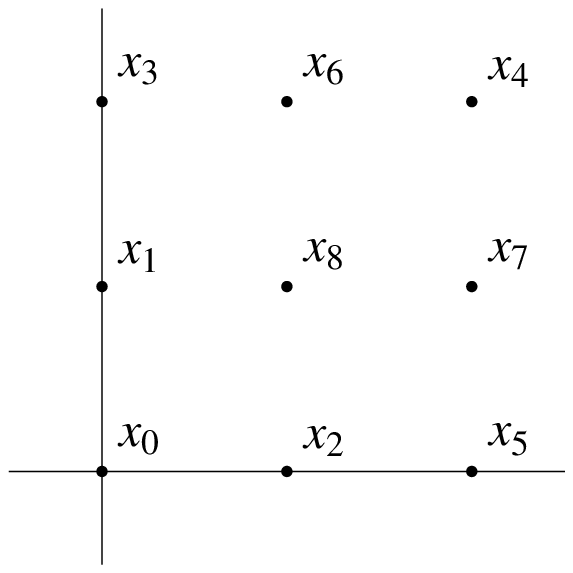}\\
The set $V(\mathbf{x})$ of values of $\Delta_{S(2)^*+(i,j)}(\mathbf{x})$ for the translated windows $S(2)^*+(i,j),\hspace{1mm}(i,j)\in T_3$ is seen to be
\begin{eqnarray}
V(\mathbf{x})=\{134, 245, 035, 467, 578, 368, 017, 128, 026\},
\end{eqnarray}
where we abbreviate the sum $x_a+x_b+x_c$ as $abc$. On the other hand the set $V(p\mathbf{x})$ of values of $\Delta_{S(2)^*+(i,j)}(\mathbf{px})$ is given by
\begin{eqnarray}
V(p\mathbf{x})=\{128, 578, 017, 368, 467, 134, 026, 245, 035\}.
\end{eqnarray}
Comparing the sets in (5.2) and (5.3) we find that $V(\mathbf{x})=V(p\mathbf{x})$. Therefore if $\mathbf{a}\in \mathbf{A}_3(\mathbb{Z}_{3^2})_{S(2)^*}$, then $p\mathbf{a}\in \mathbf{A}_3(\mathbb{Z}_{3^2})_{S(2)^*}$ too. Furthermore we can find another permutation $q$ with the same property. Let us put
\begin{eqnarray*}
q=\left(
\begin{array}{cccccccc}
 1 & 2 & 3 & 4 & 5 & 6 & 7 & 8 \\
 3 & 6 & 1 & 4 & 7 & 2 & 5 & 8
\end{array}
\right).
\end{eqnarray*}
Then the array $\mathbf{x}$ is transformed through $q$ to the following array:

\includegraphics{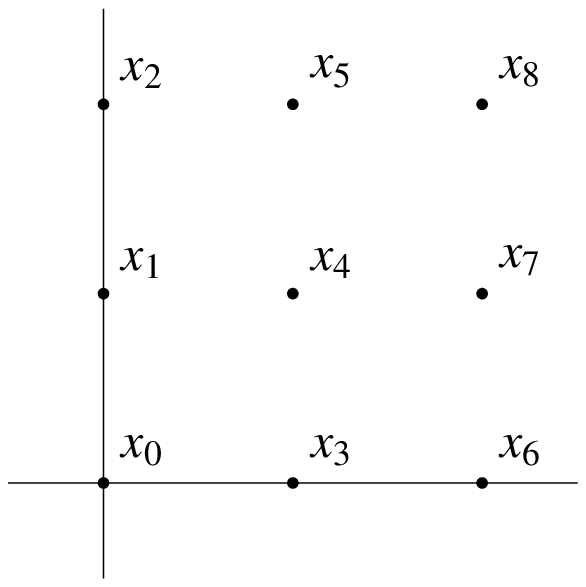}\\
Therefore the set $V(q\mathbf{x})$ of values for the transformed array $q\mathbf{x}$ is given by
\begin{eqnarray*}
V(q\mathbf{x})=\{134, 467, 017, 245, 578, 128, 035, 368, 026\}.
\end{eqnarray*}
Comparing this with (5.2),we see that $V(\mathbf{x})=V(p\mathbf{x})$. Hence if an array $\mathbf{a}$ belongs to $\mathbf{A}_3(\mathbb{Z}_{3^2})_{S(2)^*}$, then $q\mathbf{a}\in \mathbf{A}_3(\mathbb{Z}_{3^2})_{S(2)^*}$ too. Let $G(p, q)$ denote the subgroup of the symmetric group $S_8$ generated by $p$ and $q$. Since we have 
\begin{eqnarray*}
p&=&\left(
\begin{array}{cccccc}
 1 & 3 & 6 & 7 & 5 & 2 
\end{array}
\right)
\left(
\begin{array}{cc}
 4 & 8 
\end{array}
\right),\\
q&=&\left(
\begin{array}{cc}
 1 & 3
 \end{array}
\right)
\left(
\begin{array}{cc}
 2 & 6 
\end{array}
\right)
\left(
\begin{array}{cc}
 5 & 7
 \end{array}
\right),
\end{eqnarray*}
the order of $p$ (resp. $q$) is equal to six (resp. two). Furthermore we see that
\begin{eqnarray*}
qpq^{-1}&=&qpq\\
&=&\left(
\begin{array}{cccccccc}
 1 & 2 & 3 & 4 & 5 & 6 & 7 & 8 \\
 2 & 5 & 1 & 8 & 7 & 3 & 6 & 4
\end{array}
\right)\\
&=&p^{-1},
\end{eqnarray*}
Therefore the group $G(p, q)$ is isomorphic to the dihedral group $D_6$ of order 12. Thus we see that the space of balanced arrays in $\mathbf{A}_3(\mathbb{Z}_{3^2})_{S(2)^*}$ is acted upon by $G(p,q)\cong D_6$. Actually we see that the arrays in (5.1) are expressed as follows:
\begin{eqnarray*}
&&\mathbf{a}^2=qp(\mathbf{a}^1),
\mathbf{a}^3=qp^2(\mathbf{a}^1), 
\mathbf{a}^4=p(\mathbf{a}^1),\\
&&\mathbf{a}^5=qp^3(\mathbf{a}^1),
\mathbf{a}^6=p^2(\mathbf{a}^1),
\mathbf{a}^7=p^5(\mathbf{a}^1),\\
&&\mathbf{a}^8=q(\mathbf{a}^1),
\mathbf{a}^9=p^4(\mathbf{a}^1),
\mathbf{a}^{10}=qp^5(\mathbf{a}^1),\\
&&\mathbf{a}^{11}=qp^4(\mathbf{a}^1),
\mathbf{a}^{12}=p^3(\mathbf{a}^1).
\end{eqnarray*}
Thus the set of balanced arrays in $\mathbf{A}_3(\mathbb{Z}_{3^2})_{S(2)^*}$ with $x_0=0$ are acted upon transitively and faithfully by the group $G(p,q)\cong D_6$. The basic array $\mathbf{a}^1$ is illustrated as follows:\\

\noindent
$\mathbf{a}^1$:

\includegraphics{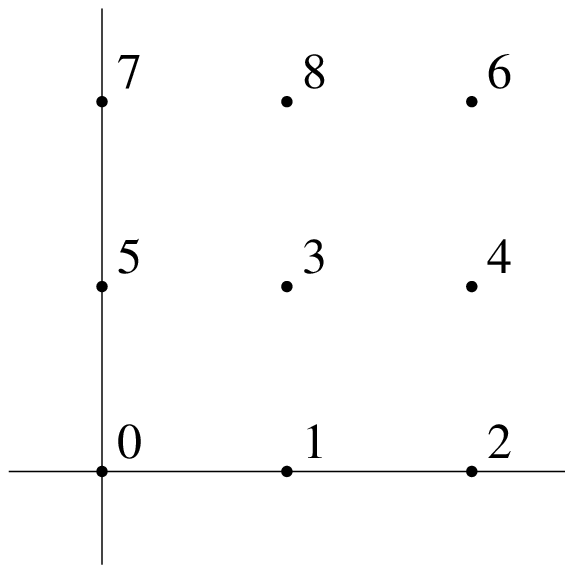}\\
We can generalize the consideration in the example above naturally, and construct a balanced array in $\mathbf{A}_n(\mathbb{Z}_{n^2})_{S(n-1)^*}$ for any $n\geq 3$:

\begin{prp}
For any $n\geq 3$, let $g_n:T_n\rightarrow\mathbb{Z}_{n^2}$ be a map defined by the rule:
\[
g_n(i,j)=
\left\{
\begin{array}{ll}
i-j, & \mbox{if}\hspace{2mm} i\geq j,\\
i-j+n, & \mbox{if}\hspace{2mm}i<j,
\end{array}
\right.
\]
and let $fn:T_n\rightarrow\mathbb{Z}_{n^2}$ be a map defined by
\begin{eqnarray*}
f_n(i,j)=g_n(i,j)+nj.
\end{eqnarray*}
Then the array $\mathbf{f}_n=(f_n(i,j))_{(i,j)\in T_n}$ is a balanced array in $\mathbf{A}_n(\mathbb{Z}_{n^2})_{S(n-1)^*}$.
\end{prp}

\begin{rem}
The array $\mathbf{a}^1$ illustrated just above {\rm Proposition 5.1} coincides with the array $\mathbf{f}_3$. The figures below depict the arrays $\mathbf{f}_4$ and $\mathbf{f}_5${\rm :}\\
$\mathbf{f}_4$:\\

\includegraphics{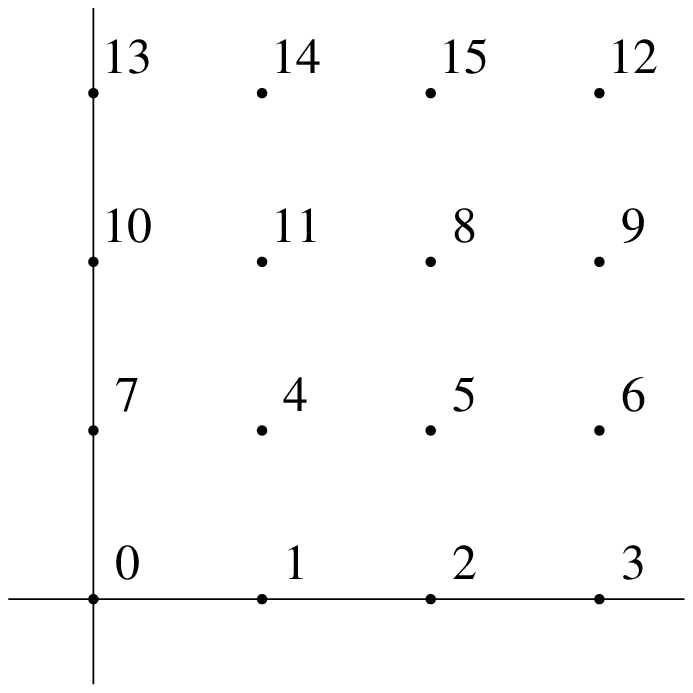}\\
$\mathbf{f}_5$:\\

\includegraphics{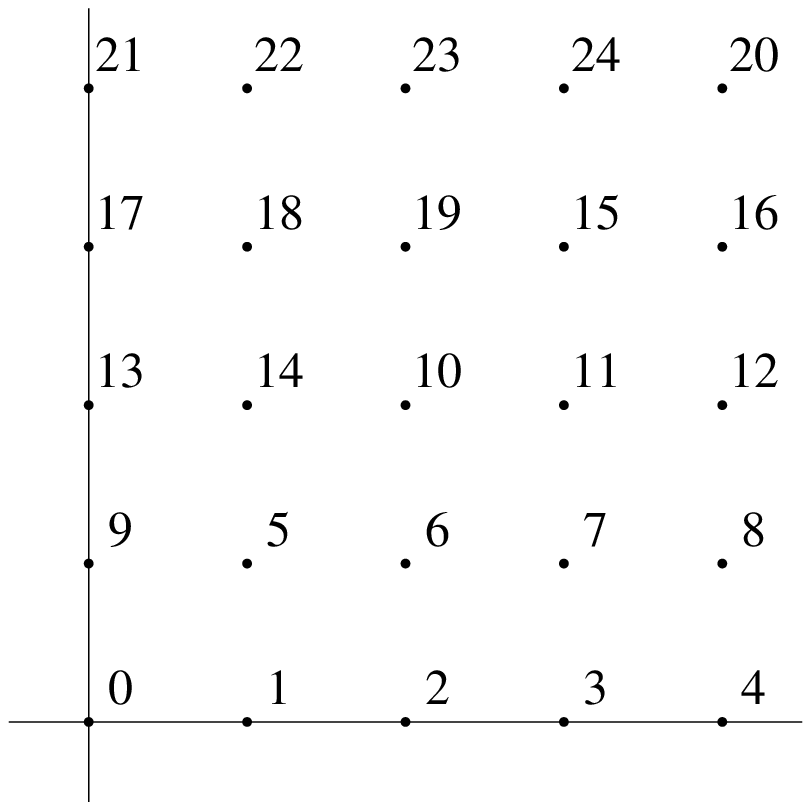}\\
One can see that the array $\mathbf{f}_n$ provides us with a natural generalization of the array $\mathbf{a}^1$.
\end{rem}

\noindent
{\it Proof of Proposition} 5.1. First we show that $\mathbf{f}_n$ is balanced, even if one might think it is evident from the figures above. Actually the equality (5.4) below which arises in the course of its proof will help us to show $\mathbf{f}_n$ is a zero-sum array. It follows from the definition of $g_n$ that the values $g_n(i,j_0)$ for a fixed $j_0\in [0, n-1]$ are as follows:
\begin{eqnarray*}
&&(g_n(i,j_0);0\leq i\leq n-1)\\
&&=(n-j_0,n-j_0+1,\cdots,n-1,0,1,\cdots,n-j_0-1).\\
\end{eqnarray*}
Therefore we have
\begin{eqnarray}
\{g_n(i,j_0);0\leq i\leq n-1\}=[0,n-1]
\end{eqnarray}
Hence the set $\{f_n(i,j_0);0\leq i\leq n-1\}$ coincides with the set $[nj_0,nj_0+n-1]$ by the definition of $f_n$. Therefore the whole set $\{f_n(i,j);0\leq i, j\leq n-1\}$ of values is equal to $[0, n^2-1]$, and hence the array $\mathbf{f}_n$ is balanced. Next we show that $\mathbf{f}_n$ is a zero-sum array for the window $S(n-1)^*$. In order to do this we consider the complement $C_n=T_n\setminus S(n-1)'$, where $S(n-1)'$ denotes the translated window $S(n-1)^*+(1,1)$. A translated set $C_n+(i_0,j_0), \hspace{1mm}(i_0,j_0)\in T_n$, consists of a pair of a horizontal and a vertical lines and a point adjacent to the intersection of the pair. More precisely we have
\begin{eqnarray}
C_n+(i_0,j_0)=H_n(i_0,j_0)\cup V_n(i_o,j_0)\cup \{(i_0+1, j_0+1)\},
\end{eqnarray}
where
\begin{eqnarray*}
H_n(i_0,j_0)&=&\{(i,j_0):i\in[0,n-1]\},\\
V_n(i_0,j_0)&=&\{(i_0,j):j\in[0,n-1]\}.
\end{eqnarray*}
Note that the horizontal line $H_n(i_0,j_0)$ and the vertical line $V_n(i_0,j_0)$ intersects at the unique point $(i_0,j_0)$:
\begin{eqnarray}
H_n(i_0,j_0)\cap V_n(i_0,j_0)=\{(i_0,j_0)\}.
\end{eqnarray}
 We compute the sum of the values of $f_n$ along these two lines. For the line $H_n(i_0,j_0)$, the sum of the values of $g_n$ is computed with the help of (5.4) as
\begin{eqnarray*}
\sum_{(i,j)\in H_n(i_0,j_0)}g_n(i,j)&=&\sum_{i\in[0,n-1]}g_n(i,j_0)\\
&=&\sum_{i\in[0,n-1]}i\\
&=&\frac{n(n-1)}{2}.
\end{eqnarray*}
Hence we have
\begin{eqnarray}
\sum_{(i,j)\in H_n(i_0,j_0)}f_n(i,j)&=&\sum_{i\in[0,n-1]}f_n(i,j_0)\nonumber\\
&=&\sum_{i\in[0,n-1]}(g_n(i,j_0)+nj_0)\nonumber\\
&=&\frac{n(n-1)}{2}+n^2j_0\nonumber\\
&=&\frac{n(n-1)}{2}.
\end{eqnarray}
(Note that we are computing everything in $\mathbb{Z}_{n^2})$. On the other hand we can show in the same way as in (5.4) that the set of values of $g_n(i,j)$ along $V_n(i_0,j_0)$ is equal to $[0,n-1]$:
\begin{eqnarray*}
\{g_n(i_0,j);0\leq j\leq n-1\}=[0,n-1]\hspace{2mm}\mbox{{\it for any}}\hspace{1mm}i_0.
\end{eqnarray*}
Therefore we have
\begin{eqnarray}
\sum_{(i,j)\in V_n(i_0,j_0)}f_n(i,j)&=&\sum_{j\in[0,n-1]}f_n(i_0,j)\nonumber\\
&=&\sum_{j\in[0,n-1]}(g_n(i_0,j)+nj)\nonumber\\
&=&\frac{n(n-1)}{2}+\frac{n^2(n-1)}{2}.
\end{eqnarray}
Combining the equalities (5.5)---(5.8), we have
\begin{eqnarray*}
&&\sum_{(i,j)\in C_n+(i_0,j_0)}f_n(i,j)\\
&&=\sum_{(i,j)\in H_n(i_0,j_0)}f_n(i,j)+\sum_{(i,j)\in V_n(i_0,j_0)}f_n(i,j)\\
&&\hspace{20mm}-f_n(i_0,j_0)+f_n(i_0+1,j_0+1)\\
&&=\left(\frac{n(n-1)}{2}\right)+\left(\frac{n(n-1)}{2}+\frac{n^2(n-1)}{2}\right)\\
&&\hspace{20mm}-(g_n(i_0,j_0)+nj_0)+(g_n(i_0+1,j_0+1)+n(j_0+1))\\
&&=\frac{n(n-1)(n+2)}{2}+n.
\end{eqnarray*}
(Note that the last equality comes from the fact that $g_n(i_0,j_0)=g_n(i_0+1,j_0+1)$ holds for any $(i_0,j_0)\in T_n$ by the definition of $g_n$.) It follows that the sum of values of $f_n$ on the translated window $S(n-1)'+(i_0,j_0)$ is computed as follows:
\begin{eqnarray*}
\sum_{(i,j)\in S(n-1)'+(i_0,j_0)}f_n(i,j)&=&\sum_{k\in [0, n^2-1]}k-\sum_{(i,j)\in C_n+(i_0,j_0)}f_n(i,j)\\
&=&\frac{n^2(n^2-1)}{2}-\left(\frac{n(n-1)(n+2)}{2}+n\right)\\
&=&\frac{n^2(n+1)(n-2)}{2}.
\end{eqnarray*}
Since $(n+1)(n-2)$ is even for any integer $n$, we see that
\begin{eqnarray*}
\sum_{(i,j)\in S(n-1)'+(i_0,j_0)}f_n(i,j)=0.\\
\end{eqnarray*}
This completes the proof of Proposition 5.1. \qed\\\\

For the square windows $S(k),\hspace{1mm}2\leq k\leq n-2$, of smaller size, the following proposition gives a partial answer for our problem (P):

\begin{prp}
For any odd prime $p\geq 5$ and for any $k$ with $2\leq k\leq p-2$, there does not exist a balanced array in $\mathbf{A}_p(\mathbb{Z}_{p^2})_{S(k)^*}$.
\end{prp}

\noindent
{\it Proof}. Contrary to the conclusion, suppose that $\mathbf{a}=(\mathbf{a}_{\mathbf{i}})$ is a balanced array in $\mathbf{A}_p(\mathbb{Z}_{p^2})_{S(k)^*}$. Let $\mathbf{b}=(\mathbf{b}_{\mathbf{i}})$ be the $\mathbb{F}_p$-valued array defined by the rule $\mathbf{b}_{\mathbf{i}}=\mathbf{a}_{\mathbf{i}}\pmod p$. Then $\mathbf{b}$ is a nonzero array and belongs necessarily to $\mathbf{A}_p(\mathbb{F}_p)_{S(k)^*}$. This contradicts to Theorem 4.1. This concludes the proof. \qed\\\\

\noindent
\begin{rem}
The author does not know the answer for the problem {\rm (P)} when $n$ is a composite number. The problem is related to the topic in Section four, which is restricted to the case of $p$-torus, and might require a good control of the behavior of the $n$-th roots of unity under the extensions of the finite field $\mathbb{F}_p$.
\end{rem}

\noindent
{\large References}\\
$[1]$ G. Frobenius, \"{U}ber die primfaktoren der gruppendeterminante, Sitzungber. Preuss. Akad. Wiss. Berlin (1896) 1343-1382.\\
$[2]$ F. Hazama, Discrete tomography and Hodge cycles, Tohoku Math. J. $\mathbf{59}(2007), 423-440$.\\
$[3]$ F. Hazama, Discrete tomography through distribution theory, Publ. RIMS, Kyoto Univ. $\mathbf{44}(2008), 1069-1095$.\\
$[4]$ M. Nivat, On a tomographic equivalence between (0,1)-matrices, in {\it Theory is Forever}, 216-234, Lecture Notes in Comput. Sci. 3113, Springer, Berlin, 2004.
\end{document}